%Paper: q-alg/9505029
%From: chered@math.unc.edu (Ivan Cherednik)
%Date: Wed, 24 May 1995 19:55:46 +0500
%Date (revised): Wed, 24 May 1995 20:14:38 +0500
%Date (revised): Tue, 6 Jun 1995 10:31:50 +0500

% formatam.tex -- AMSTeX template file
% Version: June 10, 1992

\input amstex

%\documentstyle{amams} % input Annals of Mathematics macros.
%%%%%%%%%%%%%%%%%%%%%%%%%%%%%%%%%%%%%%%%%%%%%%%%%%%%%%%%%%%
%% amams.sty: AMSTeX Macros for Articles to be published in
%%
%% Annals of Mathematics
%%
%% Princeton University and the
%% Institute for Advanced Study
%%
%% Published by Princeton University Press
%%%%%%%%%%%%%%%%%%%%%%%%%%%%%%%%%%%%%%%%%%%%%%%%%%%%%%%%%%%

%%%%%%%%%%%%%%%%%%%%%%%%%%%%%%%%%%%%%%%%%%%%%%%%%%%%%%%%%%%
%% Variations on AMSPPT.sty written by Amy Hendrickson
%% TeXnology Inc, Brookline, MA
%% 617 738-8029, amyh@ai.mit.edu
%%%%%%%%%%%%%%%%%%%%%%%%%%%%%%%%%%%%%%%%%%%%%%%%%%%%%%%%%%%

\def\spaces{\space\space\space\space\space\space\space\space\space\space}
\def\spacess{\message{\spaces\spaces\spaces\spaces\spaces\spaces\spaces}}
\spacess
\spacess
\message{Annals of Mathematics Style: Current Version: 1.1. June 10, 1992}
\spacess
\spacess
%%%%%%%%%%%%%%%%%%%%%%%%%%%%%%%%%%%%%%%%%%%%%%%%%%%%%%%%

\catcode`\@=11

\hyphenation{acad-e-my acad-e-mies af-ter-thought anom-aly anom-alies
an-ti-deriv-a-tive an-tin-o-my an-tin-o-mies apoth-e-o-ses apoth-e-o-sis
ap-pen-dix ar-che-typ-al as-sign-a-ble as-sist-ant-ship as-ymp-tot-ic
asyn-chro-nous at-trib-uted at-trib-ut-able bank-rupt bank-rupt-cy
bi-dif-fer-en-tial blue-print busier busiest cat-a-stroph-ic
cat-a-stroph-i-cally con-gress cross-hatched data-base de-fin-i-tive
de-riv-a-tive dis-trib-ute dri-ver dri-vers eco-nom-ics econ-o-mist
elit-ist equi-vari-ant ex-quis-ite ex-tra-or-di-nary flow-chart
for-mi-da-ble forth-right friv-o-lous ge-o-des-ic ge-o-det-ic geo-met-ric
griev-ance griev-ous griev-ous-ly hexa-dec-i-mal ho-lo-no-my ho-mo-thetic
ideals idio-syn-crasy in-fin-ite-ly in-fin-i-tes-i-mal ir-rev-o-ca-ble
key-stroke lam-en-ta-ble light-weight mal-a-prop-ism man-u-script
mar-gin-al meta-bol-ic me-tab-o-lism meta-lan-guage me-trop-o-lis
met-ro-pol-i-tan mi-nut-est mol-e-cule mono-chrome mono-pole mo-nop-oly
mono-spline mo-not-o-nous mul-ti-fac-eted mul-ti-plic-able non-euclid-ean
non-iso-mor-phic non-smooth par-a-digm par-a-bol-ic pa-rab-o-loid
pa-ram-e-trize para-mount pen-ta-gon phe-nom-e-non post-script pre-am-ble
pro-ce-dur-al pro-hib-i-tive pro-hib-i-tive-ly pseu-do-dif-fer-en-tial
pseu-do-fi-nite pseu-do-nym qua-drat-ics quad-ra-ture qua-si-smooth
qua-si-sta-tion-ary qua-si-tri-an-gu-lar quin-tes-sence quin-tes-sen-tial
re-arrange-ment rec-tan-gle ret-ri-bu-tion retro-fit retro-fit-ted
right-eous right-eous-ness ro-bot ro-bot-ics sched-ul-ing se-mes-ter
semi-def-i-nite semi-ho-mo-thet-ic set-up se-vere-ly side-step sov-er-eign
spe-cious spher-oid spher-oid-al star-tling star-tling-ly
sta-tis-tics sto-chas-tic straight-est strange-ness strat-a-gem strong-hold
sum-ma-ble symp-to-matic syn-chro-nous topo-graph-i-cal tra-vers-a-ble
tra-ver-sal tra-ver-sals treach-ery turn-around un-at-tached un-err-ing-ly
white-space wide-spread wing-spread wretch-ed wretch-ed-ly Brown-ian
Eng-lish Euler-ian Feb-ru-ary Gauss-ian Grothen-dieck Hamil-ton-ian
Her-mit-ian Jan-u-ary Japan-ese Kor-te-weg Le-gendre Lip-schitz
Lip-schitz-ian Mar-kov-ian Noe-ther-ian No-vem-ber Rie-mann-ian
Schwarz-schild Sep-tem-ber Za-mo-lod-chi-kov Knizh-nik quan-tum Op-dam
Mac-do-nald Ca-lo-ge-ro Su-ther-land Mo-ser Ol-sha-net-sky  Pe-re-lo-mov
in-de-pen-dent ope-ra-tors
}

\Invalid@\nofrills
\Invalid@\usualspace
\newif\ifnofrills@
\def\nofrills@#1#2{\relaxnext@
  \DN@{\ifx\next\nofrills
    \nofrills@true\let#2\relax\DN@\nofrills{\nextii@}%
  \else
    \nofrills@false\def#2{#1}\let\next@\nextii@\fi
\next@}}
\def\usualspace@#1{\ifnofrills@\def\usualspace{#1}\fi}
\def\addto#1#2{\csname \expandafter\eat@\string#1@\endcsname
  \expandafter{\the\csname \expandafter\eat@\string#1@\endcsname#2}}
\newdimen\bigsize@
\def\big@#1#2{{\hbox{$\left#2\vcenter to#1\bigsize@{}%
  \right.\nulldelimiterspace\z@\m@th$}}}
\def\big{\big@\@ne}
\def\Big{\big@{1.5}}
\def\bigg{\big@\tw@}
\def\Bigg{\big@{2.5}}
\def\raggedcenter@{\leftskip\z@ plus.4\hsize \rightskip\leftskip
 \parfillskip\z@ \parindent\z@ \spaceskip.3333em \xspaceskip.5em
 \pretolerance9999\tolerance9999 \exhyphenpenalty\@M
 \hyphenpenalty\@M \let\\\linebreak}
\def\upperspecialchars{\def\ss{SS}\let\i=I\let\j=J\let\ae\AE\let\oe\OE
  \let\o\O\let\aa\AA\let\l\L}
\def\uppercasetext@#1{%
  {\spaceskip1.2\fontdimen2\the\font plus1.2\fontdimen3\the\font
   \upperspecialchars\uctext@#1$\m@th\aftergroup\eat@$}}
\def\uctext@#1$#2${\endash@#1-\endash@$#2$\uctext@}
\def\endash@#1-#2\endash@{%
\uppercase{#1}\if\notempty{#2}--\endash@#2\endash@\fi}
\def\runaway@#1{\DN@{#1}\ifx\envir@\next@
  \Err@{You seem to have a missing or misspelled \string\end#1 ...}%
  \let\envir@\empty\fi}
\newif\iftemp@
\def\notempty#1{TT\fi\def\test@{#1}\ifx\test@\empty\temp@false
  \else\temp@true\fi \iftemp@}

%\comment%%% remove
\font@\tensmc=cmcsc10
\font@\sevenex=cmex7
\font@\sevenit=cmti7
\font@\eightrm=cmr8 % preloaded in plain.tex
\font@\sixrm=cmr6 % preloaded in plain.tex
\font@\eighti=cmmi8     \skewchar\eighti='177 % preloaded
\font@\sixi=cmmi6       \skewchar\sixi='177   % preloaded
\font@\eightsy=cmsy8    \skewchar\eightsy='60 % preloaded
\font@\sixsy=cmsy6      \skewchar\sixsy='60   % preloaded
\font@\eightex=cmex8 %
\font@\eightbf=cmbx8 % preloaded in plain.tex
\font@\sixbf=cmbx6   % preloaded in plain.tex
\font@\eightit=cmti8 % preloaded in plain.tex
\font@\eightsl=cmsl8 % preloaded in plain.tex
\font@\eightsmc=cmcsc10
\font@\eighttt=cmtt8 % preloaded in plain.tex
%\font@\ninerm=cmr9
%\font@\ninei=cmmi9    \skewchar\ninei='177
%\font@\ninesy=cmsy9   \skewchar\ninesy='60
%\font@\nineex=cmex9
%\font@\ninebf=cmbx9
%\font@\nineit=cmti9
%\font@\ninesl=cmsl9
%\font@\ninesmc=cmcsc9
%\font@\ninemsa=msam9
%\font@\ninemsb=msbm9
%\font@\nineeufm=eufm9
%\endcomment%%%

\loadmsam
\loadmsbm
\loadeufm
\UseAMSsymbols

\def\penaltyandskip@#1#2{\relax\ifdim\lastskip<#2\relax\removelastskip
      \ifnum#1=\z@\else\penalty@#1\relax\fi\vskip#2%
  \else\ifnum#1=\z@\else\penalty@#1\relax\fi\fi}
\def\nobreak{\penalty\@M
  \ifvmode\def\penalty@{\let\penalty@\penalty\count@@@}%
  \everypar{\let\penalty@\penalty\everypar{}}\fi}
\let\penalty@\penalty

\def\block{\RIfMIfI@\nondmatherr@\block\fi
       \else\ifvmode\vskip\abovedisplayskip\noindent\fi
        $$\def\endblock{\par\egroup$$}\fi
  \vbox\bgroup\advance\hsize-2\indenti\noindent}
\def\endblock{\par\egroup}

\def\logo@{\baselineskip2pc \hbox to\hsize{\hfil\eightpoint Typeset by
 \AmSTeX}}

%%%%%%%%%%%%%%%%%%%%%%%%%%%%%%%%%%%%%%%%%%%%%%%%%%%%%%%%%%%%%%%
%% Macros for Annals of Mathematics written by Amy Hendrickson
%% TeXnology Inc, Brookline, MA
%% 617 738-8029, amyh@ai.mit.edu
%%%%%%%%%%%%%%%%%%%%%%%%%%%%%%%%%%%%%%%%%%%%%%%%%%%%%%%%%%%%%%%

%% This file includes:
%% 1) Font declarations,
%% 2) Page set up,
%% 3) Title page
%% 4) Section heads,
%% 5) Equation macros, autonumbering equations, etc.,
%% 6) Figure and Table Captions,
%% 7) End matter macros: Bibliography, Appendix, etc.,
%% 8) Footnotes,
%% 9) Theorem type environments
%% 10) Cross-referencing
%% 11) Listing
%% 12) Article and Journal Table of Contents

%%%%%%%%%%%%%%%%%%%%%%%%%%%%%%%%%%%
%% 1) Font declarations,
% Computer Modern fonts

% Small Caps
\font\elevensc=cmcsc10 scaled\magstephalf
\font\tensc=cmcsc10

\font\eightsc=cmcsc10 scaled800

\font\elevenrm=cmr10 scaled \magstephalf%!!!
\font\ninerm=cmr9
\font\eightrm=cmr8
\font\sixrm=cmr6
\font\fiverm=cmr5

\font\eleveni=cmmi10 scaled\magstephalf
\font\ninei=cmmi9
\font\eighti=cmmi8
\font\sixi=cmmi6
\font\fivei=cmmi5
\skewchar\ninei='177 \skewchar\eighti='177 \skewchar\sixi='177
\skewchar\eleveni='177

\font\elevensy=cmsy10 scaled\magstephalf
\font\ninesy=cmsy9
\font\eightsy=cmsy8
\font\sixsy=cmsy6
\font\fivesy=cmsy5
\skewchar\ninesy='60 \skewchar\eightsy='60 \skewchar\sixsy='60
\skewchar\elevensy'60

\font\eighteenbf=cmbx10 scaled\magstep3

\font\twelvebf=cmbx10 scaled \magstep1
\font\elevenbf=cmbx10 scaled \magstephalf
\font\tenbf=cmbx10
\font\ninebf=cmbx9
\font\eightbf=cmbx8
\font\sixbf=cmbx6
\font\fivebf=cmbx5

\font\elevenit=cmti10 scaled\magstephalf
\font\nineit=cmti9
\font\eightit=cmti8

% Fonts for bold math
\font\eighteenmib=cmmib10 scaled \magstep3
\font\twelvemib=cmmib10 scaled \magstep1
\font\elevenmib=cmmib10 scaled\magstephalf
\font\tenmib=cmmib10
\font\eightmib=cmmib10 scaled 800
\font\sixmib=cmmib10 scaled 600

\font\eighteensyb=cmbsy10 scaled \magstep3
\font\twelvesyb=cmbsy10 scaled \magstep1
\font\elevensyb=cmbsy10 scaled \magstephalf
\font\tensyb=cmbsy10
\font\eightsyb=cmbsy10 scaled 800
\font\sixsyb=cmbsy10 scaled 600

\font\elevenex=cmex10 scaled \magstephalf
\font\tenex=cmex10
\font\eighteenex=cmex10 scaled \magstep3

%%%%%%%%%%%%%%%%%%%%%%%%%%%%
%% Font families

\def\elevenpoint{\def\rm{\fam0\elevenrm}%
  \textfont0=\elevenrm \scriptfont0=\eightrm \scriptscriptfont0=\sixrm
  \textfont1=\eleveni \scriptfont1=\eighti \scriptscriptfont1=\sixi
  \textfont2=\elevensy \scriptfont2=\eightsy \scriptscriptfont2=\sixsy
  \textfont3=\elevenex \scriptfont3=\tenex \scriptscriptfont3=\tenex
  \def\bf{\fam\bffam\elevenbf}%
  \def\it{\fam\itfam\elevenit}%
  \textfont\bffam=\elevenbf \scriptfont\bffam=\eightbf
   \scriptscriptfont\bffam=\sixbf
\normalbaselineskip=13.95pt
  \setbox\strutbox=\hbox{\vrule height9.5pt depth4.4pt width0pt\relax}%
  \normalbaselines\rm}

\elevenpoint %%% default fonts and baselineskip

\def\ninepoint{\def\rm{\fam0\ninerm}%
  \textfont0=\ninerm \scriptfont0=\sixrm \scriptscriptfont0=\fiverm
  \textfont1=\ninei \scriptfont1=\sixi \scriptscriptfont1=\fivei
  \textfont2=\ninesy \scriptfont2=\sixsy \scriptscriptfont2=\fivesy
  \textfont3=\tenex \scriptfont3=\tenex \scriptscriptfont3=\tenex
  \def\it{\fam\itfam\nineit}%
  \textfont\itfam=\nineit
  \def\bf{\fam\bffam\ninebf}%
  \textfont\bffam=\ninebf \scriptfont\bffam=\sixbf
   \scriptscriptfont\bffam=\fivebf
\normalbaselineskip=11pt
  \setbox\strutbox=\hbox{\vrule height8pt depth3pt width0pt\relax}%
  \normalbaselines\rm}

\def\eightpoint{\def\rm{\fam0\eightrm}%
  \textfont0=\eightrm \scriptfont0=\sixrm \scriptscriptfont0=\fiverm
  \textfont1=\eighti \scriptfont1=\sixi \scriptscriptfont1=\fivei
  \textfont2=\eightsy \scriptfont2=\sixsy \scriptscriptfont2=\fivesy
  \textfont3=\tenex \scriptfont3=\tenex \scriptscriptfont3=\tenex
  \def\it{\fam\itfam\eightit}%
  \textfont\itfam=\eightit
  \def\bf{\fam\bffam\eightbf}%
  \textfont\bffam=\eightbf \scriptfont\bffam=\sixbf
   \scriptscriptfont\bffam=\fivebf
\normalbaselineskip=12pt
  \setbox\strutbox=\hbox{\vrule height8.5pt depth3.5pt width0pt\relax}%
  \normalbaselines\rm}

%%%%%%%%%%%%%%%%%%%%%%%%%%%%
%% Font families for bold math in title and section heads

\def\eighteenbold{\def\rm{\fam0\eighteenbf}%
  \textfont0=\eighteenbf \scriptfont0=\twelvebf \scriptscriptfont0=\tenbf
  \textfont1=\eighteenmib \scriptfont1=\twelvemib\scriptscriptfont1=\tenmib
  \textfont2=\eighteensyb \scriptfont2=\twelvesyb\scriptscriptfont2=\tensyb
  \textfont3=\eighteenex \scriptfont3=\tenex \scriptscriptfont3=\tenex
  \def\bf{\fam\bffam\eighteenbf}%
  \textfont\bffam=\eighteenbf \scriptfont\bffam=\twelvebf
   \scriptscriptfont\bffam=\tenbf
\normalbaselineskip=20pt
  \setbox\strutbox=\hbox{\vrule height13.5pt depth6.5pt width0pt\relax}%
\everymath {\fam0 }
\everydisplay {\fam0 }
  \normalbaselines\rm}

\def\elevenbold{\def\rm{\fam0\elevenbf}%
  \textfont0=\elevenbf \scriptfont0=\eightbf \scriptscriptfont0=\sixbf
  \textfont1=\elevenmib \scriptfont1=\eightmib \scriptscriptfont1=\sixmib
  \textfont2=\elevensyb \scriptfont2=\eightsyb \scriptscriptfont2=\sixsyb
  \textfont3=\elevenex \scriptfont3=\elevenex \scriptscriptfont3=\elevenex
  \def\bf{\fam\bffam\elevenbf}%
  \textfont\bffam=\elevenbf \scriptfont\bffam=\eightbf
   \scriptscriptfont\bffam=\sixbf
\normalbaselineskip=14pt
  \setbox\strutbox=\hbox{\vrule height10pt depth4pt width0pt\relax}%
\everymath {\fam0 }
\everydisplay {\fam0 }
  \normalbaselines\bf}

%%%%%%%%%%%%%%%%%%%%%%%%%%%%%%%%%%%%%%%%%%%%%%%%%%%%%%%%%
%% 2) Page set up
\hsize=31pc
\vsize=48pc

\parindent=22pt
\parskip=0pt

\widowpenalty=10000
\clubpenalty=10000

\topskip=12pt

\skip\footins=20pt
\dimen\footins=3in % maximum footnote height

\abovedisplayskip=6.95pt plus3.5pt minus 3pt
\belowdisplayskip=\abovedisplayskip

%% Output routine

\voffset=7pt\hoffset= .7in%7pt magstep1

\newif\iftitle%!

\def\amheadline{\iftitle%
\hbox to\hsize{\hss\currannalsline\hss}\else\line{\ifodd\pageno
\hfill\thetitle\hfill\llap{\elevenrm\folio}\else\rlap{\elevenrm\folio}
\hfill\theauthors\hfill\fi}\fi}

\headline={\amheadline}%!!!
\footline={\global\titlefalse}
%\output={\bindingoffset\plainoutput}

%%%%%%%%%%%%%%%%%%%%%%%%%%%%%%%%%%%%%%%
% 3) Title page

 %#1= Volume number, #2=year of publication
\def\annalsline#1#2{\vfill\eject
\ifodd\pageno\else % first page of article on right.
\line{\hfill}
\vfill\eject\fi
\global\titletrue
\def\currannalsline{\eightrm %Annals of Mathematics,%ANNALS
{\eightbf#1} (#2), \thepages}}

\def\titleheadline#1{\def\one{#1}\ifx\one\empty\else
\def\thetitle{{%\frenchspacing%
\let\\ \relax\eightsc\uppercase{#1}}}\fi}

\newif\ifshort

\let\shorttitle\titleheadline

\def\onpages#1#2{\def\thepages{#1--#2}}

\def\thismuchskip[#1]{\vskip#1pt}
\def\ilook{\ifx\next[ \let\go\thismuchskip\else
\let\go\relax\vskip1pt\fi\go}

\def\institution#1{\def\theinstitutions{\vbox{\baselineskip10pt
\def\and{{\eightrm and }}
\def\\{\futurelet\next\ilook}\eightsc #1}}}
\let\institutions\institution

\newwrite\auxfile

\def\startingpage#1{\def\one{#1}\ifx\one\empty\global\pageno=1\else
\global\pageno=#1\fi
\theoremcount=0 \eqcount=0 \sectioncount=0
\openin1 \jobname.aux \ifeof1
\onpages{#1}{???}
\else\closein1 \relax\input \jobname.aux
\onpages{#1}{\lastpage}
\fi\immediate\openout\auxfile=\jobname.aux
}

\def\endarticle{\ifRefsUsed\global\RefsUsedfalse%
\else\vskip21pt\theinstitutions%
\nobreak\vskip8pt
%\vbox{\thereceived\therevised}%
\fi%
\write\auxfile{\string\def\string\lastpage{\the\pageno}}}

\outer\def\bye{\endarticle\par \vfill \supereject \end}

% variation on code from amsspt.sty ==>
\def\document{\let\fontlist@\relax\let\alloclist@\relax
 \elevenpoint}%%% add for annals!!!

% <=== end of code varied from amsppt.sty

\newif\ifacks
\long\def\acknowledgements#1{\def\one{#1}\ifx\one\empty\else
\vskip-\baselineskip%
\global\ackstrue\footnote{\ \unskip}{*#1}\fi}

\def\title#1{\titleheadline{#1}
\vbox to80pt{\vfill
\baselineskip=18pt
\parindent=0pt
\overfullrule=0pt
\hyphenpenalty=10000
\everypar={\hskip\parfillskip\relax}
\hbadness=10000
\def\\ {\vskip1sp}
\eighteenbold#1\vskip1sp}}

\newif\ifauthor

\def\author#1{\vskip11pt
\hbox to\hsize{\hss\tenrm By \tensc#1\ifacks\global\acksfalse*\fi\hss}
\ifshort\else\xdef\theauthors{{\eightsc\uppercase{#1}}}\fi%
\vskip21pt\global\authortrue\everypar={\global\authorfalse\everypar={}}}

\def\twoauthors#1#2{\vskip11pt
\hbox to\hsize{\hss%
\tenrm By \tensc#1 {\tenrm and} #2\ifacks\global\acksfalse*\fi\hss}
\ifshort\else\xdef\theauthors{{\eightsc\uppercase{#1 and #2}}}\fi%
\vskip21pt
\global\authortrue\everypar={\global\authorfalse\everypar={}}}

%%%%%%%%%%%%%%%%%%%%%%%%%%%%%%%%
%% 4) Section heads, counters

\newcount\theoremcount
\newcount\sectioncount
\newcount\eqcount

\newif\ifspecialnumon

\def\eqnumber=#1 {\global\eqcount=#1 \global\advance\eqcount by-1\relax}
\def\sectionnumber=#1 {\global\sectioncount=#1
\global\advance\sectioncount by-1\relax}
\def\proclaimnumber=#1 {\global\theoremcount=#1
\global\advance\theoremcount by-1\relax}

\newif\ifsection
\newif\ifsubsection

\def\elevenboldmath#1{$#1$\egroup}
\def\mathbold{\hbox\bgroup\elevenbold\elevenboldmath}

\def\section#1{\global\theoremcount=0
\global\eqcount=0
\ifauthor\global\authorfalse\else%
\vskip18pt plus 18pt minus 6pt\fi%
{\parindent=0pt
\everypar={\hskip\parfillskip}%            !!! remove
\def\\ {\vskip1sp}\elevenpoint\bf%
\ifspecialnumon\global\specialnumonfalse$\rm\spnum$%
\gdef\sectnum{$\rm\spnum$}%
\else\interlinepenalty=10000%
\global\advance\sectioncount by1\relax\the\sectioncount%
\gdef\sectnum{\the\sectioncount}%
\fi. \hskip6pt#1%                          !!!add }} and stop here
\vrule width0pt depth12pt}
\hskip\parfillskip%\break%!
\global\sectiontrue%
\everypar={\global\sectionfalse\global\interlinepenalty=0\everypar={}}%
\ignorespaces

}

%%%%%%%%%%%%%%%%%%%%%%%%%%%%%%%%
%% 5) Equation Macros

\newif\ifspequation

\let\eqno\leqno %automatic left side equation numbers %%!!!remove l-eqno

\newif\ifineqalignno
\let\saveleqalignno\leqalignno                        %%!!!remove l-eqno
\def\leqalignno{\let\eqnu\Eeqnu\saveleqalignno}

\let\eqalignno\leqalignno

\def\sectandeqnum{%
\ifspecialnumon\global\specialnumonfalse
$\rm\spnum$\gdef\eqnum{{$\rm\spnum$}}\else\global\firstlettertrue
\global\advance\eqcount by1
\ifappend\applett\else\the\sectioncount\fi.%
\the\eqcount
\xdef\eqnum{\ifappend\applett\else\the\sectioncount\fi.\the\eqcount}\fi}

\def\eqnu{\leqno{\hbox{\elevenrm\ifspequation\else(\fi\sectandeqnum
\ifspequation\global\spequationfalse\else)\fi}}}      %!!! l-eqno

\def\Speqnu{\global\setbox\leqnobox=\hbox{\elevenrm
\ifspequation\else%
(\fi\sectandeqnum\ifspequation\global\spequationfalse\else)\fi}}

\def\Eeqnu{\hbox{\elevenrm
\ifspequation\else%
(\fi\sectandeqnum\ifspequation\global\spequationfalse\else)\fi}}

\newif\iffirstletter
\global\firstlettertrue
\def\eqletter#1{\global\specialnumontrue\iffirstletter\global\firstletterfalse
\global\advance\eqcount by1\fi
\gdef\spnum{\the\sectioncount.\the\eqcount#1}\eqnu}

%%% Split math
\newbox\leqnobox
\def\outsideeqnu#1{\global\setbox\leqnobox=\hbox{#1}}

\def\eatone#1{}

%% Vertically centers equation number.
\def\dosplit#1#2{\vskip-.5\abovedisplayskip
\setbox0=\hbox{$\let\eqno\outsideeqnu%
\let\eqnu\Speqnu\let\leqno\outsideeqnu#2$}%
\setbox1\vbox{\noindent\hskip\wd\leqnobox\ifdim\wd\leqnobox>0pt\hskip1em\fi%
$\displaystyle#1\mathstrut$\hskip0pt plus1fill\relax
\vskip1pt
\line{\hfill$\let\eqnu\eatone\let\leqno\eatone%
\displaystyle#2\mathstrut$\ifmathqed~~\qed\fi}}%
\copy1
\ifvoid\leqnobox
\else\dimen0=\ht1 \advance\dimen0 by\dp1
\vskip-\dimen0
\vbox to\dimen0{\vfill
\hbox{\unhbox\leqnobox}
\vfill}
\fi}

\everydisplay{\lookforbreak}

\long\def\lookforbreak #1$${\def\mathone{#1}
\expandafter\testforbreak\mathone\splitmath @}

\def\testforbreak#1\splitmath #2@{\def\mathtwo{#2}\ifx\mathtwo\empty%
#1$$%
\ifmathqed\vskip-\belowdisplayskip
\setbox0=\vbox{\let\eqno\relax\let\eqnu\relax$\displaystyle#1$}%
\vskip-\ht0\vskip-3.5pt\hbox to\hsize{\hfill\qed}
\vskip\ht0\vskip3.5pt\fi
\else$$\vskip-\belowdisplayskip
\vbox{\dosplit{#1}{\let\eqno\eatone
\let\splitmath\relax#2}}%
\nobreak\vskip.5\belowdisplayskip
\noindent\ignorespaces\fi}

%% Proof box to be used when proof ends with equation.

\newif\ifmathqed

%%%%%%%%%%%%%%%%%%%%%%%%%%%%%
%% \mtable, Math table to make binary table easily

%% Use:
% \mtable
% &n_1&n_2&n_3&n_4&n_5&n_6\cr
% \Delta_1&M_3&M_2&0&0&0&0\cr
% \Delta_2&0&0&M_1&M_3&0&0\cr
% \endmtable

\newcount\linenum
\newcount\colnum

%++
\def\spline{\omit&\multispan{\the\colnum}{\hrulefill}\cr}
\def\colcounter{\ifnum\linenum=1\global\advance\colnum by1\fi}

\def\everyline{\noalign{\global\advance\linenum by1\relax}%
\ifnum\linenum=2\spline\fi}

\def\mtable{\bgroup\offinterlineskip
\everycr={\everyline}\global\linenum=0
\halign\bgroup\vrule height 10pt depth 4pt width0pt
\hfill$##$\hfill\hskip6pt\ifnum\linenum>1
\vrule\fi&&\colcounter\hskip12pt\hfill$##$\hfill\hskip12pt\cr}

\def\endmtable{\crcr\egroup\egroup}

%%%%%%%%%%%%%%%%%%%%%%%%%%%%%
% Array

%% Will work in math or in text, will be in math mode inside array.
%% For each column desired supply
%% r, l, or c, for right, left, or center orientation of that column.
%% End each line with \\.

%% To use:
%  \array ccc*
%  x_s\leq a_1\\
%  a_s<x_s^s<b_s\\
%  x_s\geq a_1
%  \endarray

\def\xast{*}
\newcount\intable
\newcount\mathcol
\newcount\savemathcol
\newcount\topmathcol
\newdimen\arrayhspace
\newdimen\arrayvspace

\arrayhspace=8pt % horizontal space between columns, (half this width
                 %  will horizontally precede and follow the array)
\arrayvspace=12pt % vertical space between lines

\newif\ifdollaron

\def\mathalign#1{\def\arg{#1}\ifx\arg\xast%
\let\go\relax\else\let\go\mathalign%
\global\advance\mathcol by1 %
\global\advance\topmathcol by1 %
\expandafter\def\csname  mathcol\the\mathcol\endcsname{#1}%
\fi\go}

\def\arraypickapart#1]#2*{\if#1c \ifmmode\vcenter\else
\global\dollarontrue$\vcenter\fi\else%
\if#1t\vtop\else\if#1b\vbox\fi\fi\fi\bgroup%
\def\one{#2}}

\def\arraystrut{\vrule height .7\arrayvspace depth .3\arrayvspace width 0pt}

\def\array#1#2*{\def\firstarg{#1}%
\if\firstarg[ \def\two{#2} \expandafter\arraypickapart\two*\else%
\ifmmode\vcenter\else\vbox\fi\bgroup \def\one{#1#2}\fi%
\global\everycr={\noalign{\global\mathcol=\savemathcol\relax}}%
\def\\ {\cr}%
\global\advance\intable by1 %
\ifnum\intable=1 \global\mathcol=0 \savemathcol=0 %
\else \global\advance\mathcol by1 \savemathcol=\mathcol\fi%
\expandafter\mathalign\one*%
\mathcol=\savemathcol %
\halign\bgroup&\hskip.5\arrayhspace\arraystrut%
\global\advance\mathcol by1 \relax%
\expandafter\if\csname mathcol\the\mathcol\endcsname r\hfill\else%
\expandafter\if\csname mathcol\the\mathcol\endcsname c\hfill\fi\fi%
$\displaystyle##$%
\expandafter\if\csname mathcol\the\mathcol\endcsname r\else\hfill\fi\relax%
\hskip.5\arrayhspace\cr}

\def\endarray{\crcr\egroup\egroup%
\global\mathcol=\savemathcol %
\global\advance\intable by -1\relax%
\ifnum\intable=0 %
\ifdollaron\global\dollaronfalse $\fi
\loop\ifnum\topmathcol>0 %
\expandafter\def\csname  mathcol\the\topmathcol\endcsname{}%
\global\advance\topmathcol by-1 \repeat%
\global\everycr={}\fi%
}

\def\big#1{{\hbox{$\left#1\vbox to 10pt{}\right.\n@space$}}}
\def\Big#1{{\hbox{$\left#1\vbox to 13pt{}\right.\n@space$}}}
\def\bigg#1{{\hbox{$\left#1\vbox to 16pt{}\right.\n@space$}}}
\def\Bigg#1{{\hbox{$\left#1\vbox to 19pt{}\right.\n@space$}}}

%%%%%%%%%%%%%%%%%%%%%%%%%%%%%%%%%%%%%%%%%%%%%%%%%%%%%%%%%%%%%%%%
% 6) Figure and Table Captions.

\def\figcaption#1#2#3{\topinsert
\vskip4pt %<===topadjust to match height of ascenders on opposing page.
\vbox to#3{\vfill}\vskip1sp
\setbox0=\hbox{\eightsc Figure #1.\hskip12pt\eightpoint #2}
\ifdim\wd0>\hsize
\noindent\eightsc Figure #1.\hskip12pt\eightpoint #2
\else
\centerline{\eightsc Figure #1.\hskip12pt\eightpoint #2}
\fi
\vskip16pt
\endinsert}

\def\wfig#1#2#3{\topinsert
\vskip4pt %<===topadjust to match height of ascenders on opposing page.
\hbox to\hsize{\hss\vbox{\hrule height .25pt width #3
\hbox to #3{\vrule width .25pt height #2\hfill\vrule width .25pt height #2}
\hrule height.25pt}\hss}
\vskip1sp
\centerline{\eightsc Figure #1}
\vskip16pt
\endinsert}

\def\wfigcaption#1#2#3#4{\topinsert
\vskip4pt %<===topadjust to match height of ascenders on opposing page.
\hbox to\hsize{\hss\vbox{\hrule height .25pt width #4
\hbox to #4{\vrule width .25pt height #3\hfill\vrule width .25pt height #3}
\hrule height.25pt}\hss}
\vskip1sp
\setbox0=\hbox{\eightsc Figure #1.\hskip12pt\eightpoint\rm #2}
\ifdim\wd0>\hsize
\noindent\eightsc Figure #1.\hskip12pt\eightpoint\rm #2\else
\centerline{\eightsc Figure #1.\hskip12pt\eightpoint\rm #2}\fi
\vskip16pt
\endinsert}

\def\tabcaption#1#2{\vskip6pt
\setbox0=\hbox{\eightsc Table #1.\hskip12pt\eightpoint #2}
\ifdim\wd0>\hsize
\noindent\eightsc Table #1.\hskip12pt\eightpoint #2
\else
\centerline{\eightsc Table #1.\hskip12pt\eightpoint #2}
\fi
\vskip6pt}

\def\endinsert{\egroup\if@mid\dimen@\ht\z@\advance\dimen@\dp\z@
\advance\dimen@ 12\p@\advance\dimen@\pagetotal\ifdim\dimen@ >\pagegoal
\@midfalse\p@gefalse\fi\fi\if@mid\smallskip\box\z@\bigbreak\else
\insert\topins{\penalty 100 \splittopskip\z@skip\splitmaxdepth\maxdimen
\floatingpenalty\z@\ifp@ge\dimen@\dp\z@\vbox to\vsize {\unvbox \z@
\kern -\dimen@ }\else\box\z@\nobreak\smallskip\fi}\fi\endgroup}

\def\pagecontents{
\ifvoid\topins \else\iftitle\else
\unvbox \topins \fi\fi \dimen@ =\dp \@cclv \unvbox
\@cclv
\ifvoid\topins\else\iftitle\unvbox\topins\fi\fi
\ifvoid \footins \else \vskip \skip \footins \footnoterule
\unvbox \footins \fi \ifr@ggedbottom \kern -\dimen@ \vfil \fi}

%%%%%%%%%%%%%%%%%%%%%%%%%%%%%%%%%%%%%%%%%%%%%%%%%%%%%%%%%%%%%%%%
% 7) End Matter

\newif\ifappend

\def\appendix#1#2{\def\applett{#1}\def\two{#2}%
\global\appendtrue
\global\theoremcount=0
\global\eqcount=0
\vskip18pt plus 18pt
\vbox{\parindent=0pt
\everypar={\hskip\parfillskip}
\def\\ {\vskip1sp}\elevenbold Appendix%
\ifx\applett\empty\gdef\applett{A}\ifx\two\empty\else.\fi%
\else\ #1.\fi\hskip6pt#2\vskip12pt}%
\global\sectiontrue%
\everypar={\global\sectionfalse\everypar={}}\nobreak\ignorespaces}

\newif\ifRefsUsed
\long\def\references{\global\RefsUsedtrue\vskip21pt
\theinstitutions
\global\everypar={}\global\bibnum=0
\vskip20pt\goodbreak\bgroup
\vbox{\centerline{\eightsc References}\vskip6pt}%
\ifdim\maxbibwidth>0pt
\leftskip=\maxbibwidth%
\parindent=-\maxbibwidth%
\else
\leftskip=18pt%
\parindent=-18pt%
\fi
\ninepoint
\frenchspacing
\nobreak\ignorespaces\everypar={\amref}%
}

\def\endreferences{\vskip1sp\egroup\global\everypar={}%
\nobreak\vskip8pt\vbox{\thereceived\therevised}
}

\newcount\bibnum

\def\amref#1 {\global\advance\bibnum by1%
\immediate\write\auxfile{\string\expandafter\string\def\string\csname
\space #1croref\string\endcsname{[\the\bibnum]}}%
\leavevmode\hbox to18pt{\hbox to13.2pt{\hss[\the\bibnum]}\hfill}}

\def\bibline{\hbox to30pt{\hrulefill}\/\/}

\def\name#1{{\eightsc#1}}

\newdimen\maxbibwidth
\def\AuthorRefNames [#1] {%
\immediate\write\auxfile{\string\def\string\cite\string##1{[\string##1]}}

\def\amref{\spamref}
\setbox0=\hbox{[#1] }\global\maxbibwidth=\wd0\relax}

\def\spamref[#1] {\leavevmode\hbox to\maxbibwidth{\hss[#1]\hfill}}

%%%%%%%%%%%%%%%%%%%%%%%%%%%%%%%%%%%%%%%%%%%%%%%%%%%%%%%%%%%%%%%%
%% 8) Footnotes

\def\footnoterule{\kern-3pt\hrule width1in height.5pt\kern2.5pt}

\def\footnote#1#2{%
\plainfootnote{#1}{{\eightpoint\normalbaselineskip11pt
\normalbaselines#2}}}

\def\vfootnote#1{%
\insert \footins \bgroup \eightpoint\baselineskip11pt
\interlinepenalty \interfootnotelinepenalty
\splittopskip \ht \strutbox \splitmaxdepth \dp \strutbox \floatingpenalty
\@MM \leftskip \z@skip \rightskip \z@skip \spaceskip \z@skip
\xspaceskip \z@skip
{#1}$\,$\footstrut \futurelet \next \fo@t}

%%%%%%%%%%%%%%%%%%%%%%%%%%%%%%%%%%%%%%%%%%%%%%%%%%%%%%%%%%%%%%%%
%% 9) Theorem type environments

\newif\iffirstadded
\newif\ifadded

\def\addedlett{}

\def\alltheoremnums{%
\ifspecialnumon\global\specialnumonfalse
\ifadded\global\addedfalse
\iffirstadded\global\firstaddedfalse
\global\advance\theoremcount by1 \fi
\ifappend\applett\else\the\sectioncount\fi.\the\theoremcount\addedlett%
\xdef\theoremnum{\ifappend\applett\else\the\sectioncount\fi.%
\the\theoremcount\addedlett}%
\else$\rm\spnum$\def\theoremnum{{$\rm\spnum$}}\fi%
\else\global\firstaddedtrue
\global\advance\theoremcount by1
\ifappend\applett\else\the\sectioncount\fi.\the\theoremcount%
\xdef\theoremnum{\ifappend\applett\else\the\sectioncount\fi.%
\the\theoremcount}\fi}

\def\allcorolnums{%
\ifspecialnumon\global\specialnumonfalse
\ifadded\global\addedfalse
\iffirstadded\global\firstaddedfalse
\global\advance\corolcount by1 \fi
\the\corolcount\addedlett%
\else$\rm\spnum$\def\corolnum{$\rm\spnum$}\fi%
\else\global\advance\corolcount by1
\the\corolcount\fi}

%% use for Theorem, Corollary, Lemma, Proposition, Demonstration and similar.

\newcount\corolcount
\def\xcorol{Corollary}
\def\xtheorem{Theorem}
\def\xmaintheorem{Main Theorem}

\newif\ifthtitle

\let\saverparen)
\let\savelparen(
\def\rmparenl{{\rm(}}
\def\rmparenr{{\rm\/)}}
{
\catcode`(=13
\catcode`)=13
\gdef\makeparensRM{\catcode`(=13\catcode`)=13\let(=\rmparenl%
\let)=\rmparenr%
\everymath{\let(\savelparen%
\let)\saverparen}%
\everydisplay{\let(\savelparen%
\let)\saverparen\lookforbreak}}}

\medskipamount=8pt plus.1\baselineskip minus.05\baselineskip

\def\rmtext#1{\hbox{\rm#1}}

\def\proclaim#1{\vskip-\lastskip
\def\one{#1}\ifx\one\xtheorem\global\corolcount=0\fi
\ifsection\global\sectionfalse\vskip-6pt\fi
\medskip
{\elevensc#1}%
\ifx\one\xmaintheorem\global\corolcount=0
\gdef\theoremnum{Main Theorem}\else%
\ifx\one\xcorol\ \allcorolnums\else\ \alltheoremnums\fi\fi%
\ifthtitle\ \global\thtitlefalse{\rm(\thethtitle)}\fi.%
\hskip1em\bgroup\let\text\rmtext\makeparensRM\it\ignorespaces}

\def\nonumproclaim#1{\vskip-\lastskip
\def\one{#1}\ifx\one\xtheorem\global\corolcount=0\fi
\ifsection\global\sectionfalse\vskip-6pt\fi
\medskip
{\elevensc#1}.\ifx\one\xmaintheorem\global\corolcount=0
\gdef\theoremnum{Main Theorem}\fi\hskip.5pc%
\bgroup\it\makeparensRM\ignorespaces}

\def\endproclaim{\egroup\medskip}

%% Use demo for Proof, Proof of, Definition, Example,
%% Remark, Case, Subcase, Conjecture, Note, Notation,
%% Convention, Construction and Step.
%% Any other use for demo will format similar to `Proof.'

\def\xproof{Proof}
\def\xremark{Remark}
\def\xcase{Case}
\def\xsubcase{Subcase}
\def\xconjecture{Conjecture}
\def\xstep{Step}
\def\xof{of}

\def\deconstruct#1 #2 #3 #4 #5 @{\def\one{#1}\def\two{#2}\def\three{#3}%
\def\four{#4}%
\ifx\two\empty #1\else%
\ifx\one\xproof%
\ifx\two\xof%
  \ifx\three\xcorol Proof of Corollary \rm#4\else%
     \ifx\three\xtheorem Proof of Theorem \rm#4\else\xone\fi%
  \fi\fi%
\else\xone\fi\fi.}

\def\pickup#1 {\def\this{#1}%
\ifx\this\xproof\global\let\go\demoproof
\global\let\enddemo\endproof\else
\ifx\this\xremark\global\let\go\demoremark\else
\ifx\this\xcase\global\let\go\demostep\else
\ifx\this\xsubcase\global\let\go\demostep\else
\ifx\this\xconjecture\global\let\go\demostep\else
\ifx\this\xstep\global\let\go\demostep\else
\global\let\go\demoproof\fi\fi\fi\fi\fi\fi}

\newif\ifnonum
\def\demo#1{\vskip-\lastskip
\ifsection\global\sectionfalse\vskip-6pt\fi
\def\one{#1 }\def\two{#1*}%
\setbox0=\hbox{\expandafter\pickup\one}\expandafter\go\two}

\def\numbereddemo#1{\vskip-\lastskip
\ifsection\global\sectionfalse\vskip-6pt\fi
\def\two{#1*}%
\expandafter\demoremark\two}

\def\demoproof#1*{\medskip\def\xone{#1}
{\ignorespaces\it\expandafter\deconstruct\xone {} {} {} {} {} @%
\unskip\hskip6pt}\rm\ignorespaces}

\def\demoremark#1*{\medskip
{\it\ignorespaces#1\/} \ifnonum\global\nonumtrue\else
 \alltheoremnums\unskip.\fi\hskip1pc\rm\ignorespaces}

\def\demostep#1 #2*{\vskip4pt
{\it\ignorespaces#1\/} #2.\hskip1pc\rm\ignorespaces}

\def\enddemo{\medskip}

\def\endproof{\ifmathqed\global\mathqedfalse\medskip\else
\parfillskip=0pt~~\hfill\qed\medskip
\fi\global\parfillskip0pt plus 1fil\relax
\gdef\enddemo{\medskip}}

\def\qed{\vbox{\hrule\hbox{\vrule height6pt\hskip6pt\vrule}\hrule}}

%% Proof box to be used in a \proclaim{}...\endproclaim environment

\def\proofbox{\parfillskip=0pt~~\hfill\qed\vskip1sp\parfillskip=
0pt plus 1fil\relax}

%%%%

%%%%%%%%%%%%%%%%%%%%%
%% 10) CrossRefs

%%% Generic crossreferencing
%%% to use: \label\nameoflabel* (will give the page number when referenced)

% Commands to access current state of counter, for cross-referencing
% \sectnum
% \theoremnum
% \eqnum

%%% You can make another definition that includes counters and/or the
%%% page number and access this information as the second argument:
%%% \label\yourlabelname[2.13]*

%%% Since this method of cross-referencing relies
%%% on an auxiliary file, the first time you tex the file
%%% you will get `??' when you write \ref\nameoflabel.
%%% When you TeX the file the second time the auxiliary file
%%% will be input and \ref\nameoflabel will produce the cross-ref.

\def\stripbs#1#2*{\def\one{#2}}

\def\emptyspace{ }
\def\nextthing{}
\def\newline{***}
\def\eatone#1{ }

\def\lookatspace#1{\ifcat\noexpand#1\ \else%
\gdef\nextthing{}\xdef\next{#1}%
\ifx\next\emptyspace%
\let\nextthing\emptyspace\else\ifx\next\newline%
\gdef\nextthing{\eatone}\fi\fi\fi\egroup\nextthing#1}

{\catcode`\^^M=\active%
\gdef\spacer{\bgroup\catcode`\^^M=\active%
\let^^M=\newline\obeyspaces\lookatspace}}

\def\ref#1{\seeifdefined{#1}\expandafter\csname\one\endcsname\spacer}

\def\cite#1{\expandafter\ifx\csname#1croref\endcsname\relax[??]\else
\csname#1croref\endcsname\fi\spacer}

%% for testing in \label and \ref to see if term already labeled.

\def\seeifdefined#1{\expandafter\stripbs\string#1croref*%
\crorefdefining{#1}}

\newif\ifcromessage
\global\cromessagetrue

\def\crorefdefining#1{\ifdefined{\one}{}
{\ifcromessage\global\cromessagefalse%
\message{\spaces\spaces\spaces\spaces\spaces\spaces\spaces}%
\message{<Undefined reference.}%
\message{Please TeX file once more to have accurate cross-references.>}%
\message{\spaces\spaces\spaces\spaces\spaces\spaces\spaces}\fi[??]}}

\def\label#1#2*{\gdef\ctest{#2}%
\xdef\currlabel{\string#1croref}
\expandafter\seeifdefined{#1}%
\ifx\empty\ctest%
\xdef\labelnow{\write\auxfile{\noexpand\def\currlabel{\the\pageno}}}%
\else\xdef\labelnow{\write\auxfile{\noexpand\def\currlabel{#2}}}\fi%
\labelnow}

\def\ifdefined#1#2#3{\expandafter\ifx\csname#1\endcsname\relax%
#3\else#2\fi}

%%%%%%%%%%%%%%%%%%%%%
%% 11) Listing

%% To use with asterisks:

%%%%%%%%%%%%%%%%%%%%%%
%% 12) Article and Journal Table of Contents

\def\articlecontents{
\vskip20pt\centerline{\bf Table of Contents}\everypar={}\vskip6pt
\bgroup \leftskip=3pc \parindent=-2pc
\def\item##1{\vskip1sp\indent\hbox to2pc{##1.\hfill}}}

\def\endcontents{\vskip1sp\leftskip=0pt\egroup}

\def\journalcontents{\vfill\eject
\def\currannalsline{\hfill}
\global\titletrue
\vglue3.5pc
\centerline{\tensc\hskip12pt TABLE OF CONTENTS}\everypar={}\vskip30pt
\bgroup \leftskip=34pt \rightskip=-12pt \parindent=-22pt
  \def\\ {\vskip1sp\noindent}
\def\pagenum##1{\unskip\parfillskip=0pt\dotfill##1\vskip1sp
\parfillskip=0pt plus 1fil\relax}
\def\name##1{{\tensc##1}}}

%% default values

\institution{}
\onpages{0}{0}
\def\lastpage{???}
\def\thetitle{Title ???}
\def\theauthors{Authors ???}
\def\thereceived{}
\def\therevised{}

\gdef\split{\relaxnext@\ifinany@\let\next\insplit@\else
 \ifmmode\ifinner\def\next{\onlydmatherr@\split}\else
 \let\next\outsplit@\fi\else
 \def\next{\onlydmatherr@\split}\fi\fi\let\eqnu\xspliteqnu\next}

\gdef\align{\relaxnext@\ifingather@\let\next\galign@\else
 \ifmmode\ifinner\def\next{\onlydmatherr@\align}\else
 \let\next\align@\fi\else
 \def\next{\onlydmatherr@\align}\fi\fi\let\eqnu\xspliteqnu\next}

\def\spliteqnu{{\tenrm\sectandeqnum}\relax}

\def\xspliteqnu{\tag\spliteqnu}

\catcode`@=12

\document

%-------------- Publisher's entries --------------------
\annalsline{May}{1995}
%\line{\hfil Revised }
\startingpage{1}     %numeration
%%\received{??}
%\revised{February, 1995}

%\magnification=\magstep1

%--------------- Author macros ---------------
%                   MACROS
%
%                                 AUX
%
%
%
%                      endaux
%

\def\iif{\quad\hbox{ if }\quad}

\def\for{\  \hbox{ for } \ }
\def\if{ \ \hbox{ if } \ }
\def\when{ \ \hbox{ when } \ }
\def\where{\  \hbox{ where } \ }
\def\and{\  \hbox{ and } \ }

\def\equal{\buildrel  def \over =}

\def\la{\lambda}
\def\La{\Lambda}
\def\om{\omega}
\def\Om{\Omega}

\def\th{\theta}
\def\al{\alpha}
\def\be{\beta}
\def\ga{\gamma}
\def\ep{\epsilon}

\def\de{\delta}

\def\si{\sigma}
\def\Si{\Sigma}
\def\Ga{\Gamma}
\def\ze{\zeta}

    %from copy, ell
\def\pa{\partial}

\def\vph{\varphi}

\def\vep{\varepsilon}

\def\tal{\tilde{\alpha}}
\def\tbe{\tilde{\beta}}
\def\tep{\tilde{\epsilon}}
\def\tde{\tilde{\delta}}
\def\tpi{\tilde{\pi}}
\def\tPi{\tilde{\Pi}}
\def\tV{\tilde{V}}

\def\tw{\tilde w}

\def\tB{\tilde B}

\def\tV{\tilde V}
\def\tz{\tilde z}
\def\tb{\tilde b}
\def\ta{\tilde a}

\def\hH{\hat{H}}

\def\hY{\hat{Y}}

\def\hT{\hat{T}}

\def\hw{\hat{w}}
\def\hu{\hat{u}}

\def\hv{\hat{v}}

\def\C{\bold{C}}
\def\Q{\bold{Q}}
\def\B{\bold{B}}
\def\R{\bold{R}}
\def\N{\bold{N}}
\def\Z{\bold{Z}}

\def\one{\bold{1}}

\def\0{\bold{0}}

\def\C{\hbox{\bf C}}

\def\H{\bold{H}}% macdonald

\def\f{\Cal{F}}
\def\t{\Cal{T}}

\def\l{\Cal{L}}

\def\p{\Cal{P}}

\def\h{\Cal{H}}

\def\v{\Cal{V}}

\font\germ=eufb10 %at 12pt
%\font\germm=eufb9 at 12pt
%\font\germ=eufm9 at 12pt
\def\goth#1{\hbox{\germ #1}}

\def\TT{\goth{T}}
\def\HH{\goth{H}}

\def\BB{\goth{B}}
\def\AA{\goth{A}}

\font\smm=msbm10 at 12pt
\def\symbol#1{\hbox{\smm #1}}
\def\lsmash{{\symbol n}}
\def\rsmash{{\symbol o}}

%endmacros

%------------------------------------------------------------------
%-------------- Author entries --------------------

%\comment                               %to remove the title
\title
{Non-symmetric Macdonald's polynomials}

 %Article title
\shorttitle{ Non-symmetric polynomials}
 % Shortened version for headline title

% Acknowledgements: Please enter all acknowledgements here.
\acknowledgements{
Partially supported by NSF grant DMS--9301114}

% Please uncomment and use appropriate command:
\author{ Ivan Cherednik}
%\twoauthors{}{}
%\authors{}% Separate each author with a comma and a space.

% Institution:
\institutions{
Math. Dept, University of North Carolina at Chapel Hill,
 N.C. 27599-3250
\\ Internet: chered\@math.unc.edu
}

%\endcomment                              %to remove the title
%-------------- Article Text--------------------
%\intro %(Optional, Introduction)
%
%
%
%                        INTRO
%
%
%{\bf 0. Introduction.
\vfil
Recently Eric Opdam [O2] in the differential case and then Ian
Macdonald [M3] in the  difference $q,t$-setting introduced remarkable
orthogonal polynomials. In contrast to major known families   they
linearly generate the space of all (non-symmetric) polynomials. Their
meaning still needs to be clarified. At the present time we believe in
their importance mainly because they are eigenfunctions of the
differential  [C4,C5] and difference [C1,C2] Dunkl operators. The
latter operators play a preponderant role in the representation theory
of (affine and) double affine Hecke algebras and related harmonic
analysis.  Anyway the Hecke algebra technique works  better for
non-symmetric polynomilas than for their celebrated symmetric
counterparts defined in [M1,M2].

Paper [O2] is mostly about analytic theory of graded affine Hecke
algebras defined by Lusztig. Its first algebraic part is  a
generalization (and a simplification) of [O1] where the Macdonald
conjectures were proved in the differential case. In particular, it
contains the presentation of the symmetric (Jacobi) polynomials in
terms of non-symmetric ones, the formula for the norms of the latter,
and the interpretation of the shift operator [O1,He] via the
anti-symmetric polynomials.  Macdonald announced in [M3] the difference
analogues of these results (when $t=q^k, k\in \Z$).

In the present paper, we prove  Macdonald's statements for arbitrary
$q,t$, also establishing the duality-evaluation theorem, the recurrence
theorem, and the basic facts on non-symmetric polynomials at roots of
unity, including a description of the projective action of $SL_2(\Z)$
on them (generalizing the representations from [Ki]).  We follow [C3]
devoted to the same questions in the symmetric case. The main point is
the definition of the difference spherical Fourier transform based on
the double affine Hecke algebras (generalizing the Hankel transform
from [D,J]). For instance, it readily results in the norm-formulas
 conjectured by Macdonald and proved in [C2]  and clarifies why  they
are so surprisingly simple.  Hopefully,  this version of  the Macdonald
theory can be extended to the elliptic case (see [C6]) without serious
difficulties.

Once the Fourier transform appeared, we cannot
 restrict ourselves to symmetric functions anymore.  Even  the
classical multi-component Fourier transform requires  at least the
coordinate functions and the corresponding differentiations.  It
reveals itself at many levels.

First, it is  easier to operate with the double affine Hecke algebra
than with its (very complicated) subalgebra of symmetric operators.
Second, promising applications are expected in arithmetic, where the
symmetric elliptic functions have no particular importance.  Although
much was done by means of the characters of Kac-Moody algebras (see
[K]), certainly they and their $q,t$-analogues are  not enough.  Then,
non-symmetric polynomials seem more relevant to incorporate the
Ramanujan $_1\Psi_1$-summation and its generalizations into the
Macdonald theory.  As to physics, they  can be transformed into
eigenfunctions of the so-called spin-Calogero-Sutherland hamiltonians
[C5] and its difference counterparts.  We also mention [O2,HO], which
contain a lot of analitic evidence on the same point.

In spite of all these, there should exist   deeper relations to the
representation theory and the combinatorics.  Till now, there hasn't
been any interpretation of the non-symmetric Opdam-Macdonald
polynomials as characters or generalized chearacters [EK] (the
equivalence of the spin-CS model and the affine Knizhnik-Zamolodchikov
equations from [C5] indicates that it could exists). Our technical
achievements  are far ahead of the understanding of their true place.

The author thanks G. Heckman, D. Kazhdan and A. Kirillov, Jr., and E.
Opdam for  usuful discussion.  The paper was started at UC at San
Diego.  I am grateful to A. Garsia and my colleagues for the kind
invitation and hospitality.

%
%
%		Section 1
%
%
%\vskip 10pt
\section { Affine root systems}
Let $R=\{\al\}   \subset \R^n$ be a root system of type $A,B,...,F,G$
with respect to a euclidean form $(z,z')$ on $\R^n \ni z,z'$,
$W$ the Weyl group  generated by  the reflections $s_\al$.
We assume that $(\al,\al)=2$ for long $\al$.
Let us  fix the set $R_{+}$ of positive  roots ($R_-=-R_+$),
the corresponding simple
roots $\al_1,...,\al_n$, and  their dual counterparts
$a_1 ,..., a_n,  a_i =\al_i^\vee, \where \al^\vee =2\al/(\al,\al)$.
The dual fundamental weights
$b_1,...,b_n$  are determined from the relations  $ (b_i,\al_j)=
\de_i^j $ for the
Kronecker delta. We will also introduce the dual root system
$R^\vee =\{\al^\vee, \al\in R\}, R^\vee_+$, and the lattices
$$
\eqalignno{
& A=\oplus^n_{i=1}\Z a_i \subset B=\oplus^n_{i=1}\Z b_i,
}
$$
  $A_\pm, B_\pm$  for $\Z_{\pm}=\{m\in\Z, \pm m\ge 0\}$
instead of $\Z$. (In the standard notations, $A= Q^\vee,\
B = P^\vee $ - see [B].)  Later on,
$$
\eqalign{
&\nu_{\al}=\nu_{\al^\vee}\ =\ (\al,\al),\  \nu_i\ =\ \nu_{\al_i}, \
\nu_R\ = \{\nu_{\al}, \al\in R\},
}
$$
$$\eqalign{
&\rho_\nu\ =\ (1/2)\sum_{\nu_{\al}=\nu} \al \ =
\ (\nu/2)\sum_{\nu_i=\nu}  b_i, \for\al\in R_+,\cr
&r_\nu\ =\ \rho_\nu^\vee \ =\ (2/\nu)\rho_\nu\ =\
\sum_{\nu_i=\nu}  b_i,\quad 2/\nu=1,2,3.
}
\eqnu
\label\rhor\eqnum*
$$

The vectors $\ \tal=[\al,k] \in
\R^n\times \R \subset \R^{n+1}$
for $\al \in R, k \in \Z $
form the {\it affine root system}
$R^a \supset R$ ( $z\in \R^n$ are identified with $ [z,0]$).
We add  $\al_0 \equal [-\th,1]$ to the  simple roots
for the {\it maximal root} $\th \in R$.
The corresponding set $R^a_+$ of positive roots coincides
with $R_+\cup \{[\al,k],\  \al\in R, \  k > 0\}$.

We denote the Dynkin diagram and its affine completion with
$\{\al_j,0 \le j \le n\}$ as the vertices by $\Ga$ and $\Ga^a$.
Let $m_{ij}=2,3,4,6$\  if $\al_i\and\al_j$ are joined by 0,1,2,3 laces
respectively.
The set of
the indices of the images of $\al_0$ by all
the automorphisms of $\Ga^a$ will be denoted by $O$ ($O=\{0\}
\for E_8,F_4,G_2$). Let $O^*={r\in O, r\neq 0}$.
The elements $b_r$ for $r\in O^*$ are the so-called minuscule
weights ($(b_r,\al)\le 1$ for
$\al \in R_+$).

Given $\tal=[\al,k]\in R^a,  \ b \in B$, let
$$
\eqalignno{
&s_{\tal}(\tz)\ =\  \tz-(z,\al^\vee)\tal,\
\ b'(\tz)\ =\ [z,\ze-(z,b)]
&\eqnu
%&(1.1)
}
$$
for $\tz=[z,\ze] \in \R^{n+1}$.

The {\it affine Weyl group} $W^a$ is generated by all $s_{\tal}$
(we write $W^a = <s_{\tal}, \tal\in R_+^a>)$. One can take
the simple reflections $s_j=s_{\al_j}, 0 \le j \le n,$ as its
generators and introduce the corresponding notion of the
length. This group is
the semi-direct product $W\lsmash A'$ of
its subgroups $W=<s_\al,
\al \in R_+>$ and $A'=\{a', a\in A\}$, where
$$
\eqalignno{
& a'=\ s_{\al}s_{[\al,1]}=\ s_{[-\al,1]}s_{\al}\for a=\al^{\vee},
\ \al\in R.
&\eqnu
%&(1.2)
}
$$

The {\it extended Weyl group} $ W^b$ generated by $W\and B'$
(instead of $A'$) is isomorphic to $W\lsmash B'$:
$$
\eqalignno{
&(wb')([z,\ze])\ =\ [w(z),\ze-(z,b)] \for w\in W, b\in B.
&\eqnu
}
$$

 Given $b_+\in B_+$, let
$$
\eqalignno{
&\om_{b_+} = w_0w^+_0  \in  W,\ \pi_{b_+} =
b'_+(\om_{b_+})^{-1}
\ \in \ W^b, \ \om_i=\om_{b_i},\pi_i=\pi_{b_i},
&\eqnu
\label\wo\eqnum*
}
$$
where $w_0$ (respectively, $w^+_0$) is the longest element in $W$
(respectively, in $ W_{b_+}$ generated by $s_i$ preserving $b_+$)
relative to the
set of generators $\{s_i\}$ for $i >0$.

To describe $W^b$ as an extension of $W^a$ we need the
elements $\pi_r=\pi_{b_r}, r \in O$. They leave $\Ga^a$ invariant
and form a group denoted by $\Pi$,
 which is isomorphic to $B/A$ by the natural
projection $\{b_r \to \pi_r\}$. As to $\{\om_r\}$,
they preserve the set $\{-\th,\al_i, i>0\}$.
The relations $\pi_r(\al_0)= \al_r= (\om_r)^{-1}(-\th)
$ distinguish the
indices $r \in O^*$. Moreover (see e.g. [C2]):
$$
\eqalignno{
& W^b  = \Pi \lsmash W^a, \where
  \pi_rs_i\pi_r^{-1}  =  s_j \if \pi_r(\al_i)=\al_j,\  0\le j\le n.
&\eqnu
}
$$

We extend the notion
of the length to $W^b$.
Given $\nu\in\nu_R,\  r\in O^*,\  \tw \in W^a$, and a reduced
decomposition $\tw\ =\ s_{j_l}...s_{j_2} s_{j_1} $ with respect to
$\{s_j, 0\le j\le n\}$, we call $l\ =\ l(\hw)$ the {\it length} of
$\hw = \pi_r\tw \in W^b$. Setting
$$
\eqalign{
\la(\hw) = &\{ \tal^1=\al_{j_1},\
\tal^2=s_{j_1}(\al_{j_2}),\
\tal^3=s_{j_1}s_{j_2}(\al_{j_3}),\ldots \cr
&\ldots,\tal^l=\tw^{-1}s_{j_l}(\al_{j_l}) \},
}
\eqnu
$$
\label\tal\eqnum*
one can represent
$$
\eqalign
{
&l=|\la(\hw)|=\sum_\nu l_\nu, \for l_\nu = l_\nu(\hw)=|\la_\nu(\hw)|,\cr
&\la_\nu(\hw) = \{\tal^{m},\ \nu(\tal^{m})= \nu(\tal_{j_m})= \nu\},
1\le m\le l,
}
\eqnu
\label\laset\eqnum*
$$
where $|\ |$  denotes the  number of elements,
 $\nu([\al,k]) \equal \nu_{\al}$.

To interpret the length geometrically, let us introduce
the following (affine)
action of $W^b$ on $z \in \R^n$:
$$
\eqalign{
& (wb')\langle z \rangle \ =\ w(b+z),\ w\in W, b\in B,\cr
& s_{\tal}\langle z\rangle\ =\ z - ((z,\al)+k)\al^\vee,
\ \tal=[\al,k]\in R^a,}
 \eqnu
\label\afaction\eqnum*
$$
and the affine Weyl chamber:
$$
\eqalignno{
&C^a\ =\ \bigcap_{j=0}^n L_{\al_j},\ L_{\tal}=\{z\in \R^n,\
(z,\al)+k>0 \}.
}
$$
Then (see e.g. [C2]):
$$
\eqalign{
\la_\nu(\hw)\ & =\ \{\tal\in R^a, \ \langle  C^a \rangle
\not\subset \hw\langle L_{\tal}\rangle, \ \nu({\tal})=\nu \} \cr
& =\ \{\tal\in R^a, \ l_\nu( \hw s_{\tal}) < l_\nu(\hw) \}.
}
 \eqnu
\label\lambda\eqnum*
$$
It coincides with (\ref\laset) due to the relations
$$
\eqalign{
&  \la_\nu(\hw\hu) = \la_\nu(\hu) \cup
\hu^{-1}(\la_\nu(\hw)),\
\la_\nu(\hw^{-1}) = -\hw\langle \la_\nu(\hw)\rangle\cr
&\if l_\nu(\hw\hu)=
l_\nu(\hw)+l_\nu(\hu).
}
\eqnu
\label\ltutw\eqnum*
$$

 Let us generalize the definition of $\om_{b_+},
 \pi_{b_+}$. See [C2], Definition 1.1, Proposition 1.3, and Theorem 1.4.

\proclaim{Proposition} Given $ b\in B$, the decomposition $b= \pi_b\om_b,
\om_b \in W$
can be uniquely determined  from the following equivalent conditions

{i)\ \  }
$\om_b(b) \ = \ b_-\in B_-$ and
$(\al,b) \neq 0 \if (-\al) \in R_+ \ni \om_b(\al)$,

{ii)\ }
$\om_b(b) \in B_-$ and $l(\om_b)$ is the smallest possible,

{iii) }
$l(\pi_b)+l(\om_b)\ =\ l(b)$ and $l(\om_b)$ is the biggest possible,

{iv)\ }
$\la(\om_b)\ =\ \{\al\in R_+,\ (\al,b)>0\}.$
\endproclaim
\label\PIOM\theoremnum*
\proofbox

It results from i) and iii), that
$\om_b\pi_b=b_-$. As to (\ref\wo), $\om_{b_+}$ really
sends $b_+$ to $B_-$
and is the smallest with this property.
We will also apply the formulas (see ibid.):
$$
\eqalignno{
&l_\nu(b)\ =\ \sum_{\al} |(b,\al)|,\  \al\in R_+,
\nu_\al=\nu \in \nu_R, &\eqnu
\cr
\label\lb\eqnum*
  &\la(b) = \{ \tal, \al \in R_+,
( b, \al ) > k\ge 0\} \cup
\{ \tal, \al \in  R_-, ( b, \al ) \ge k>0 \},
&\eqnu
\label\lambb\eqnum*  \cr
&\la(\pi_b) = \{ \tal,\ \al\in R_-,
&\eqnu
\label\lambpi\eqnum* \cr
&( b, \al )>k> 0 \if (\al,b)<0,\  ( b, \al )\ge k > 0 \if (\al,b)>0 \},
}
$$
where $\tal=[\al,k]\in R^a_+, b\in B$,
 $|\ | = $ is the
absolute value.

{\it Convexity.} Let us introduce two orderings
on $B$. Here and further $b_{\pm}$ are the unique elements
from $B_{\pm}$ which belong to the orbit $W(b)$. Namely,
$b_-=\om_b\pi_b=\om_b(b)$, $b_+=w_0(b_-)= \om_{-b}(b).$
So the equality   $c_-=b_- $ (or $c_+=b_+ $) means that $b,c$
belong to the same orbit.
Set
$$
\eqalignno{
&b \le c, c\ge b \for b, c\in B \iif c-b \in A_+,
&\eqnu
\label\order\eqnum*
\cr
&b \preceq c, c\succeq b \iif b_-\le c_- \hbox{\ \ or\ \ }
b_-=c_- \hbox{\ and\ } b\le c.
&\eqnu
\label\succ\eqnum*
}
$$
We  use $<,>,\prec, \succ$ respectively if $b \neq c$.
For instance,
$$c\succ b_+\Leftrightarrow b_+>c>b_-,\ \
c\succeq b_-\Leftrightarrow c\in W(b_-) \hbox{\ or\ }
c\succ b_+.
$$
The following sets
$$
\eqalign{
&\si(b)\equal \{c\in B, c\succeq b\},\
\si_*(b)\equal \{c\in B, c\succ b\}, \cr
&\si_+(b)\equal \{c\in B, c_->b_-\}\ =\ \si_*(b_+).
}
\eqnu
\label\cones\eqnum*
$$
are convex. Moreover $\si_+$ is $W$-invariant.
By convex, we mean that if
$ c, d= c+r\al^{\vee}\in \si$
for $\al\in R_+, r\in \Z_+$, then
$$
\eqalignno{
&\{c,\ c+\al^{\vee},...,c+(r-1)\al^{\vee},\ d\}\subset \si.
&\eqnu
\label\convex\eqnum*
}
$$
Actually all the elements from $\si(b)$
strictly between $c$ and $d$ (i.e.
$c+q\al,\ $ $0<q<r$) belong to $\si_+(b)$.
Let us adapt Proposition 1.5 from [C2]
to the setup of the present paper:

\proclaim{Proposition }
{i)  }
Given $ b \in B,\ \tal=[\al,k]\in \la(b)$, let $b=\hw \langle 0
\rangle,\ \hw= b s_{\tal},\ c = \hw\langle 0\rangle $. Then
$\hw=cw$ for $w\in W$ and  $c \in \si(b)$. More exactly,
$$
\eqalign{
& \{\al>0,k>0\}\  \Rightarrow\  b>c>s_\al(b)\  \Rightarrow \ c\in \si_+(b)\cr
& \{\al<0,k<(\al,b)\}\ \Rightarrow\ s_\al(b)>c>b
\Rightarrow \ c\in \si_+(b)\cr
& \{\al>0,k=0\}\  \Rightarrow\  c=b,\  \{\al<0,k=(\al,b)\}\
\Rightarrow \ c= s_\al(b)>b.
}
\eqnu
\label\bstal\eqnum*
$$

{ii)}
Let $\hw= b s_{\tal^{i_1}}...s_{\tal^{i_m}}$ ,where we
take $\tal^i$ from (\ref\tal) for
any sequence $1\le i_1<i_2<\ldots< i_m\le l=l(b)$.
Then $c=\hw\langle 0\rangle \in \si(b)$.  Moreover, $c\in \si_*(b)$
iff at least one $\tal^{i_p}=[\al,k]$ has $k>0$, and $c\in \si_+(b)$
iff at least one of them has $0<k<(b,a)$.
\label\BSTAL\theoremnum*
\endproclaim
\proofbox

The odering $\succ$ on $b\in B$
and the proposition (in a bit different but equivalent form)
were applied in  [C2] to discribe the structure of
the operators $Y_b$ (see below). On the other hand, this
odering appeared in [O2] (and then recently in [M3]) to introduce
the non-symmetric orthogonal polynomials. The coincidence
of these orderings is not by chance. It results from the
duality and the Recurrence Theorem below.

%
%
%
%                  Section 2
%
%
%\vskip 10pt
\section{ Double affine Hecke algebras}
We put
$m=2 \for D_{2k} \and C_{2k+1},\ m=1 \for C_{2k}, B_{k}$,
otherwise $m=|\Pi|$. Let us set
$$
\eqalignno{
&   t_{\tal} = t_{\nu(\tal)},\ t_j = t_{\al_j},
\where \tal \in R^a, 0\le j\le n, \cr
& X_{\tb}\ =\ \prod_{i=1}^nX_i^{k_i} q^{ k}
\if \tb=[b,k],
&\eqnu \cr
&\for b=\sum_{i=1}^nk_i b_i\in B,\ k \in {1\over m}\Z.
}
$$
\label\Xde\eqnum*
Here  $ q,\{ t_\nu , \nu \in \nu_R \},$ $X_1,\ldots,X_n$
are considered as  independent variables.

 Later on $ \C_{ q}$
 is the field of rational
functions in $ q^{1/m},$
$\C_q[X] = \C_q[X_b]$  means the algebra of
polynomials in terms of $X_i^{\pm 1}$
with the coefficients depending
on $ q^{1/m}$ rationally. We replace $ \C_{ q}$
by $ \C_{ q,t}$ if the functions (coefficients)
also depend rationally
on $\{t_\nu^{1/2} \}$.

Let $([a,k],[b,l])=(a,b)$ for $a,b\in B,\
[\al,k]^\vee=2[\al,k]/(\al,\al),\
 a_0=\al_0,\ \nu_{\al^\vee}=\nu_\al$.
We also introduce the map
$O^*\ni r\to r^*,\ \al_{r^*} \equal \pi_r^{-1}(\al_0)$.

\proclaim{Definition }
 The  double  affine Hecke algebra $\HH\ $
(see [C1,C2])
is generated over the field $ \C_{ q,t}$ by
the elements $\{ T_j,\ 0\le j\le n\}$,
pairwise commutative $\{X_b, \ b\in B\}$ satisfying (\ref\Xde),
 and the group $\Pi$ where the following relations are imposed:

(o)\ \  $ (T_j-t_j^{1/2})(T_j+t_j^{-1/2})\ =\ 0,\ 0\ \le\ j\ \le\ n$;

(i)\ \ \ $ T_iT_jT_i...\ =\ T_jT_iT_j...,\ m_{ij}$ factors on each side;

(ii)\ \   $ \pi_rT_i\pi_r^{-1}\ =\ T_j \if \pi_r(\al_i)=\al_j$;

(iii)\  $T_iX_b T_i\ =\ X_b X_{a_i}^{-1} \if (b,\al_i)=1,\
1 \le i\le  n$;

(iv)\  $T_0X_b T_0\ =\ X_{s_0(b)}\ =\ X_b X_{\th} q^{-1}
\if (b,\th)=-1$;

(v)\ \ $T_iX_b\ =\ X_b T_i$ if $(b,\al_i)=0 \for 0 \le i\le  n$;

(vi)\ $\pi_rX_b \pi_r^{-1}\ =\ X_{\pi_r(b)}\ =\ X_{\om^{-1}_r(b)}
 q^{(b_{r^*},b)},\  r\in O^*$.
\label\double\theoremnum*
\endproclaim
\proofbox

Given $\tw \in W^a, r\in O,\ $ the product
$$
\eqalignno{
&T_{\pi_r\tw}\equal \pi_r\prod_{k=1}^l T_{i_k},\where
\tw=\prod_{k=1}^l s_{i_k},
l=l(\tw),
&\eqnu
\label\Tw\eqnum*
}
$$
does not depend on the choice of the reduced decomposition
(because $\{T\}$ satisfy the same ``braid'' relations as $\{s\}$ do).
Moreover,
$$
\eqalignno{
&T_{\hv}T_{\hw}\ =\ T_{\hv\hw}\  \hbox{ whenever}\
 l(\hv\hw)=l(\hv)+l(\hw) \for
\hv,\hw \in W^b.
 %&(2.7)}
&\eqnu}
$$
\label\TT\eqnum*
  In particular, we arrive at the pairwise
commutative elements
$$
\eqalignno{
& Y_{b}\ =\  \prod_{i=1}^nY_i^{k_i} \if
b=\sum_{i=1}^nk_ib_i\in B,\where
 Y_i\equal T_{b_i},
&\eqnu
\label\Yb\eqnum*
}
$$
satisfying the relations
$$
\eqalign{
&T^{-1}_iY_b T^{-1}_i\ =\ Y_b Y_{a_i}^{-1} \if (b,\al_i)=1,
\cr
& T_iY_b\ =\ Y_b T_i \if (b,\al_i)=0, \ 1 \le i\le  n.}
%\eqno(2.9)
\eqnu
$$
Let us introduce the following elements from
$\C_t^n$:
$$
\eqalign{
&t^{\pm\rho}\equal (l_t(b_1)^{\pm 1},\ldots,l_t(b_n)^{\pm 1}),\where\cr
&l_t(\hw)\equal \ \prod_{\nu\in\nu_R} t_\nu^{l_\nu(\hw)/2},\
\hw\in W^b,
}
\eqnu
\label\qlen\eqnum*
$$
and the corresponding {\it evaluation maps}:
$$
\eqalign{
&X_i(t^{\pm\rho})= l_t(b_i)^{\pm 1} = Y_i(t^{\pm\rho}),\ 1\le i\le n.
}
\eqnu
\label\eval\eqnum*
$$
For instance, $X_{a_i}(t^{\rho})\ =\ l_t(a_i)= t_i$ (see (\ref\lb)).

We will establish the  duality of non-symmetric polynomials applying
the following theorem ([C2],[C3]).
\proclaim {Theorem}
i) The elements $H \in \HH\ $  have
the unique decompositions
$$
\eqalignno{
&H =\sum_{w\in W }  g_{w}  T_{w} f_w,\
g_{w} \in \C_{ q,t}[X],\ f_{w} \in \C_{ q,t}[Y].
&\eqnu
}
$$

ii) The   map
$$
\eqalign{
 \vph: &X_i \to Y_i^{-1},\ \  Y_i \to X_i^{-1},\  \ T_i \to T_i, \cr
&t_\nu \to t_\nu,\
 q\to  q,\ \nu\in \nu_R,\ 1\le i\le n.
}
\eqnu
\label\vph\eqnum*
$$
can be extended to an anti-involution
($\vph(AB)=\vph(B)\vph(A)$)
 of \HH\ .

iii) The  linear functional on \HH\
$$
\eqalignno{
&[\![ \sum_{w\in W }  g_{w}  T_{w} f_{w}]\!]\ =\
\sum_{w\in W} g_{w}(t^{-\rho}) l_t(w) f_{w}(t^{\rho})
&\eqnu
\label\brack\eqnum*
}
$$
is invariant with respect to $\vph$. The bilinear form
$$
\eqalignno{
&[\![ G,H]\!]\equal [\![ \vph(G)H]\!],\
G,H\in \HH\ ,
&\eqnu
\label\form\eqnum*
}
$$
is symmetric ($[\![ G,H]\!]= [\![ H,G]\!]$)
and non-degenerate.
\endproclaim
\label\dual\theoremnum*
\proofbox

The map $\vph$ is the composition of the
involution (see [C1])
$$
\eqalign{
  \vep:\ &X_i \to Y_i,\ \  Y_i \to X_i,\  \ T_i \to T_i^{-1}, \cr
&t_\nu \to t_\nu^{-1},\
 q\to  q^{-1},\ 1\le i\le n,
}
\eqnu
\label\vep\eqnum*
$$
and the main anti-involution  from [C2]
$$
\eqalign{
  & X_i^*\ =\  X_i^{-1},\   Y_i^*\ =\  Y_i^{-1},\
 T_i^* \ =\  T_i^{-1}, \cr
&t_\nu \to t_\nu^{-1},\
 q\to  q^{-1},\ 0\le i\le n.
}
\eqnu
\label\star\eqnum*
$$

Let us give the explicit formulas  for the action of $\vph,\vep$
on $T_0$:
$$
\eqalign{
&\vph(T_0)\ =\ Y_\th^{-1}T_0X_\th^{-1}\ =\ T_{s_\th}^{-1}X_\th^{-1},\cr
&\vep(T_0)\ =\ X_\th T_0^{-1}Y_\th\ =\ X_\th T_{s_\th}.
}
\eqnu
\label\vpht0\eqnum*
$$

The next theorem from [C3] will be used to obtain a projective
 action
of $GL_2(\Z)$ on the restricted non-symmetric polynomials when
$q,t$ are roots of unity.

\proclaim{ Theorem}
 {i)} Adding   $ q^{1/2m}$, the following maps
can be uniquely extended to automorphisms of \HH\ , preserving each of
$T_1,\ldots,T_n,t$ and $ q$:
$$
\eqalign{
& \tau_+: \ X_b \to X_b,\ \ Y_r \to X_rY_r q^{-(b_r,b_r)/2},\
Y_\th \to X_0^{-1}T_0^{-2}Y_\th,      \cr
& \tau_-: \ Y_b \to Y_b,\ \ X_r \to Y_r X_r q^{(b_r,b_r)/2},\
X_\th \to T_0X_0Y_\th^{-1} T_0, \cr
& \om: Y_b \to X_b^{-1},\ X_r \to X_r^{-1}Y_r X_r q^{(b_r,b_r)},\
\ X_\th \to T_0^{-1}Y_\th^{-1}T_0, \cr
&\where b\in B,\ r\in O^*,\ X_0\ =   q X_\th^{-1}.
}
\eqnu
\label\tauom\eqnum*
$$
Moreover,
$$
\eqalign{
&\tau_-\  =\  \vep\tau_+\vep =  \vph\tau_+\vph,\ \ \om\  =\
\tau_+^{-1}\tau_-\tau_+^{-1} =
\tau_-\tau_+^{-1}\tau_-.
}
\eqnu
\label\taumin\eqnum*
$$

ii) The above maps give   automorphisms of \HH\  and
the  (elliptic
braid) group \BB\  generated by
the elements $\{X_b,Y_b,T_i,\pi_r, q^{1/2m}\}$ satisfying
the relations (i)-(vi) from Definition \ref\double  and
(\ref\Yb). Let  $\AA_o$ be the group  of its automorphisms modulo
the conjugations by the elements
from  the center $Z(\B)$ of the group
$\B$ generated by $\{T_1,\ldots,T_n\}$.
Considering the images of $\vep$ (see (\ref\vep)),$\tau_{\pm},\om$
in $\AA_o$ we obtain the homomorphism $GL_2(\Z)\to \AA_o$:
$$
\eqalign{
& \Bigl(\matrix  0 &-1\\ -1& 0\endmatrix \Bigr) \to \vep,\
  \Bigl(\matrix 1& 1\\ 0& 1 \endmatrix \Bigr) \to \tau_+,\cr
& \Bigl(\matrix 0& -1\\ 1& 0\endmatrix \Bigr) \to \om,\
\Bigl(\matrix 1& 0\\ 1& 1   \endmatrix    \Bigr) \to \tau_-.
}
\eqnu
\label\glz\eqnum*
$$
\label\GLZ\theoremnum*
\endproclaim

%
%
%		Section 3
%
%
%\vskip 10pt
\section { Basic representation}
 Setting
$$
\eqalignno{
& x_{\tb}=  \prod_{i=1}^nx_i^{k_i} q^{ k} \if
\tb=[b,k],
b=\sum_{i=1}^nk_i b_i\in B,\ k \in {1\over m}\Z,
&\eqnu
\label\xde\eqnum*}
$$
for independent $x_1,\ldots,x_n$, we
 consider $\{X\}$ as  operators acting in $\C_{q,t}[x]=$
$\C_{q,t}[x_1^{\pm 1},$ $\ldots,x_n^{\pm 1}]$:
$$
\eqalignno{
& X_{\tb}(p(x))\ =\ x_{\tb} p(x),\    p(x) \in
\C_{q,t} [x].
&\eqnu}
$$
\label\X\eqnum*
The elements $\hw \in W^b$ act in $\C_{ q}[x]$
 by the
formulas:
$$
\eqalignno{
&\hw(x_{\tb})\ =\ x_{\hw(\tb)}.
&\eqnu}
$$
 In particular:
$$
\eqalignno{
&\pi_r(x_{b})\ =\ x_{\om^{-1}_r(b)} q^{(b_{r^*},b)}
\for \al_{r^*}\ =\ \pi_r^{-1}(\al_0), \ r\in O^*.
&\eqnu}
$$
\label\pi\eqnum*

The {\it Demazure-Lusztig operators} (see
[KL1, KK, C1], and [C2] for more detail )
$$
\eqalignno{
&\hT_j\  = \  t_j ^{1/2} s_j\ +\
(t_j^{1/2}-t_j^{-1/2})(X_{a_j}-1)^{-1}(s_j-1),
\ 0\le j\le n.
&\eqnu
\label\Demaz\eqnum*
}
$$
act   in $\C_{ q,t}[x]$ naturally.
We note that only $\hT_0$ depends on $ q$:
$$
\eqalign{
&\hT_0\  =  t_0^{1/2}s_0\ +\ (t_0^{1/2}-t_0^{-1/2})
( q X_{\th}^{-1} -1)^{-1}(s_0-1),\cr
&\where
s_0(X_i)\ =\ X_iX_{\th}^{-(b_i,\th)} q^{(b_i,\th)}.
}
%\eqno(2.12)
\eqnu
$$

\proclaim{Theorem }
 The map $ T_j\to \hT_j,\ X_b \to X_b$ (see (\ref\Xde,\ref\X)),
$\pi_r\to \pi_r$  (see (\ref\pi)) induces a $ \C_{ q,t}$-linear
homomorphism from \HH\ to the algebra of linear endomorphisms
of $\C_{ q,t}[x]$.
 This representation is faithful and
remains faithful when   $  q,t$ take  any non-zero
values assuming that
 $ q$ is not a root of unity (see [C2]). The image $\hat{H}$
is uniquely determined from the following condition:
$$
\eqalign{
&\hat{H}(f(x))\ =\ g(x)\for H\in \HH\ ,\if Hf(X)-g(X)\ \in\cr
 & \{\sum_{i=0}^n H_i(T_i-t_i)+
\sum_{r\in O^*} H_r(\pi_r-1), \where H_i,H_r\in \HH\ \}.
}
\eqnu
\label\hat\eqnum*
$$
\endproclaim
\proofbox
\label\faith\theoremnum*

We will also reformulate
 Lemmas 2.4 and 3.2 from
[C2] in the next proposition.
 They result directly from
Proposition \ref\BSTAL.

\proclaim{Proposition}
i) Given $b\in B$,
$$
\eqalignno{
&Y_b  =  \pi_b\ga_b^{\#b}T_{\om_b} +
\sum_{c\succ b, w\in W}  (cw) g_{b}^{cw},
\where g_{b}^{cw} \in\C_{q,t}(X),
&\eqnu
\label\ybg\eqnum*
\cr
&\ga_b^{\#b} =
\prod _{\tal \in \la(\pi_b)}
{t_{\tal}^{1/2}X_{\tal^\vee}^{-1} - t_{\tal}^{-1/2}\over
X_{\tal^\vee}^{-1} - 1} =
\prod_{a\in R_+^\vee}
{t_{a}^{1/2}(q_a)^k X_{a}^{-1} - t_{\tal}^{-1/2}\over
(q_a)^k X_{a}^{-1}  - 1}\for \cr
&(a^\vee,b)>0\Rightarrow (a^\vee,b_-)<k<0,\
(a^\vee,b)<0\Rightarrow (a^\vee,b_-)\le k<0.
}
$$

ii) Given $b\in B$,
$$
\eqalignno{
&\hY_{-b}\  =\  T_{\om_b}^{-1}\phi_b^{\#b}\pi_b^{-1} +
\sum_{c\succ b, w\in W} f_{b}^{cw} (cw)^{-1}\for
f_{b}^{cw} \in\C_{q,t}(X),
 &\eqnu
\label\ybf\eqnum*
\cr
&\phi_b^{\#b}(X;t)\ =\ \ga_b^{\#b}(X;t^*)\ =\
\ga_b^{\#b}(X^{-1};t),\ t_\nu^*= t_\nu^{-1}.
}
$$
\endproclaim
\label\YBGF\theoremnum*
 {\it Proof.} We will remind that the main step is the
presentation:
$$
\eqalignno{
&\hY_b = bG^{\bullet}_{\ta^l}\cdots G_{\ta^1}^{\bullet},\
\ta^1=\al_{j_1}^\vee,
\ta^2=s_{j_1}(\al_{j_2}^\vee),
\ta^3=s_{j_1}s_{j_2}(\al_{j_3}^\vee),\ldots,
&\eqnu
\label\ybprod\eqnum*
}
$$
where
$b=\pi_r s_{j_l}\cdots s_{j_1}$, $\ l=l(b), r\in O,
\ \ta=\tal^\vee=\tal/\nu_\al$,
$$
\eqalignno{
&G_{\ta;t} \ =\ G_{\ta}\ =\ t_{\ta}^{1/2}+
(t_{\ta}^{1/2}-t_{\ta}^{-1/2})
(X_{\ta}^{-1}-1)^{-1}(1-s_{\ta}),
&\eqnu
\label\G\eqnum*
\cr
&G^{\bullet}_{\ta}\ =\ G_{\ta} \if \al>0,\
 G^{\bullet}_{\ta}\ =\ G_{-\ta}^{-1}\ =\
G_{-\ta;t^*}\hbox{\ otherwise\ }.
}
$$
\proofbox

We note that (\ref\ybprod) also provides that
 the coefficents of the operator $\hY_b$
of any $\hw\neq b$ are zero
at the point $ \diamondsuit \equal ( X_1=...=X_n=0 )$.
Here the order of the coefficents (from $\C_{q,t}(X)$)
 and $\hw$ does not
matter since $\diamondsuit$ is $W^b$-invariant.
Thus
$$
\eqalign{
\hY_b(\diamondsuit)\ =\
\prod_\nu t_\nu^{(b,\rho_\nu)} b.
}
\eqnu
\label\diamond\eqnum*
$$
Next, we need an extended version of Proposition 3.6
from [C2].

\proclaim{Proposition }
The operators
$\{ Y_i, 1\le i\le n\}$
 preserve $\Si(b)\equal
\oplus_{c\in \si(b)}\C_{q,t} x_{c}$ and
 the $\Si_*(b)$ (defined for $\si_*(b)$)
for arbitrary $b\in B$. The operators
$\{T_j, 0\le j\le n\}$ preseve $\Si_+(b)=\Si_*(b_+)$:
$$
\eqalign{
\hT_j(x_b) &\hbox{\ mod\ } \Si_+(b)\  \cr
&=\ t_j^{1/2}s_j(x_b)+(t_j^{1/2}-t_j^{-1/2})x_b
 \if (b,\al_j)<0,\cr
&=\ t_j^{-1/2}s_j(x_b) \if (b,\al_j)>0,
\ =\ t_j^{1/2}x_b\if (b,\al_j)=0.
}
\eqnu
\label\tonx\eqnum*
$$
\endproclaim
\label\TONX\theoremnum*
{\it Proof.} Formulas (\ref\tonx) are verified
directly. Similarly, assuming that $\al>0$ in
$\tal=[\al,k]$,
$$
\eqalign{
G_{\tal^\vee}(x_b) \hbox{\ mod\ } \Si_*(b)\ &
= t_\al^{1/2}x_b+(t_\al^{1/2}-t_\al^{-1/2})s_{\tal}(x_b)
 \if (b,\al)\le 0,\cr
 &= t_\al^{-1/2}x_b \if (b,\al)>0.
}
\eqnu
\label\gonx\eqnum*
$$
Replacing  $\tal, \ t_\al$ by
 $-\tal,\ t_\al^{-1}$ we can use the same formulas
for $G^{\bullet}_{\tal^\vee}$ when $\al<0$.
\proofbox

 Relations (\ref\tonx) and
the formulas $\pi_r(x_b)=x_{\pi_r(b)}$
for $r\in O$ induce
an action of the {\it affine Hecke algebra} $\h_Y$
generated by $T_j (0\le j\le n)$ and the group $ \Pi$ in the
space
$$
V(b_-)\equal\Si(b_-)/\Si_+(b_-)\ =\ \C_{q,t}[W(b)]\ =\
\C_{q,t}\otimes \C[W(b)].
\eqnu
\label\vbmin\eqnum*
$$

The $\h_Y$-module $V(b_-)$ is irreducible (for generic $q$).
It can be decribed as the induced  representaion
generated by $x_+(=x_{b_+})$ and
satisfying the following defining conditions:
$$
\eqalign{
&T_j(x_+)\ =\ t_j^{1/2}\if (b,\al_j)=0,\ 1\le j\le n,\cr
&Y_a(x_+)\ =\ q^{(a,b_+)}\prod_\nu t_\nu^{(\om_{b_+}(a),\rho_\nu)}x_+,
\ a\in B.
}
\eqnu
\label\xplus\eqnum*
$$

We also note that  (\ref\tonx) can be rewritten in
$V(b_-)$ as follows:
$$
\eqalign{
s_j(x_b)\ &=\ (t_j^{1/2}T_j(x_b))^{-1} \if (b,\al_j)<0,\cr
&=\ t_j^{1/2}T_j(x_b) \if (b,\al_j)>0.
}
\eqnu
\label\sonx\eqnum*
$$
Hence  this representation    corresponds
to the natural action
of $W^b$ on the indices of $T_{\hw}$ after
proper normalization.

%
%
%		Section 4
%
%
%\vskip 10pt
\section { Orthogonal polynomials}
The coefficient of $x^0=1$ ({\it the constant term})
of a polynomilal
$f\in  \C_{q,t}[x]$
will be denoted by $\langle  f \rangle$. Let
$$
\eqalign{
&\mu\ =\ \prod_{a \in R_+^\vee}
\prod_{i=0}^\infty {(1-x_aq_a^{i}) (1-x_a^{-1}q_a^{i+1})
\over
(1-x_a t_aq_a^{i}) (1-x_a^{-1}t_a^{}q_a^{i+1})},
}
\eqnu
\label\mu\eqnum*
$$
where $q_a=q_{\nu}=q^{2/\nu} \for \nu=\nu_a$.

The coefficients of $\mu_1\equal \mu/\langle \mu \rangle$
are from $\C(q,t)$, where the formula for the
constant term of $\mu$ is as follows
(see [C2]):
$$
\eqalign{
&\langle\mu\rangle\ =\ \prod_{a \in R_+^\vee}
\prod_{i=1}^\infty {(1-x_a(t^\rho)q_a^{i})^2
\over
(1-x_a(t^\rho) t_aq_a^{i}) (1-x_a(t^\rho) t_a^{-1}q_a^{i})}.
}
\eqnu
\label\consterm\eqnum*
$$
Here
 $x_{b}(t^{\pm\rho}q^{c})=
q^{(b,c)}\prod_\nu t_\nu^{\pm(b,\rho_\nu)}$.

We note that
$\mu_1^*\ =\ \mu_1$  with respect to the involution
$$
 x_b^*\ =\  x_{-b},\ t^*\ =\ t^{-1},\ q^*\ =\ q^{-1}.
$$
If $t_\nu=q_\nu^{k_\nu}$ for $k_\nu\in \Z_+$ then $\mu\in
\C(q,t)[x]$ (see (\ref\resinner)).

Setting
$$
\eqalignno{
&\langle f,g\rangle\ =\langle \mu_1 f\ {g}^*\rangle\ =\
\langle g,f\rangle^* \for
f,g \in \C(q,t)[x],
&\eqnu
\label\innerpro\eqnum*
}
$$
 we  introduce the {\it non-symmetric Macdonald
polynomials} $e_b(x),\   b \in B_-$, by means of
the conditions
$$
\eqalignno{
&e_b-x_b\ \in\ \Si_*(b),\
\langle e_b, x_{c}\rangle = 0 \for c\in \si_* =
\{c\in B, c\succ b\}
&\eqnu
\label\macd\eqnum*
}
$$
in the setup of Section 1.
They can be determined by the Gram - Schmidt process
because the  pairing
is non-degenerate
 and form a
basis in $\C(q,t)[x]$. We also note that $w_0(e_b(x^{-1}))=e_{w_0(-b)}$
since  $-w_0$ does not change
the odering $\succ$.

This definition is due to Macdonald [M3] who generalized
 Opdam's non-symmetric polynomials introduced
in the degenerate (differential) case in [O2]. He also
established the connection with the $Y$-operators from
the previous section, which will be  discussed next.
 In Opdam's paper,
the trigonometric Dunkl operators from [C5] play the role
 of
$\{Y_b\}$.

 The notations are
from  Proposition \ref\PIOM and (\ref\rhor).
 We
identify the operators $H\in \HH\ $ with their images
$\hH$ and use the involution
$ \bar{x}_a=x_a^{-1},\
\bar{q}= q,\ \bar{t}=t,\ a\in B$.

\proclaim {Theorem}
The polynomials $\{e_b,b\in B\}$ are eigenvectors of
 the operators $\{L_f\equal f(Y_1,\cdots, Y_n), f\in \C_{q,t}[x]\}$:
$$
\eqalignno{
&L_{\bar{f}}(e_b)\ =\ f(\#b)e_b, \where
\#b\equal\pi_b=b\om_b^{-1},
&\eqnu
\label\Yone\eqnum*
\cr
& x_a(bw)\equal x_a(q^b t^{-w(\rho)})\ =\
q^{(a,b)}\prod_\nu t_\nu^{-(w(\rho_\nu),a)},\ w\in W.
}
$$
\endproclaim
\label\YONE\theoremnum*
{\it Proof.} Due to [C2],
$\langle Hf,g\rangle  = $ $ \langle f,H^*g\rangle$
for any $H\in \HH\ $ for the anti-involution $^*$
from (\ref\star). Hence
the operators $\{Y_b\}$ are unitary
relative to $\langle\ ,\ \rangle$.
Since they leave all
$\Si(a),\Si_*(a)$ invariant
(Proposition \ref\TONX), their eigenvectors in $\C_{q,t}[x]$ are
exactly $\{e\}$. The eigenvalues are readily calculated
by means of formula (\ref\diamond).
\proofbox

It is worth mentioning that
$$
\eqalignno{
&x_a(\#b)\ =\
q^{(a,b)}\prod_\nu t_\nu^{-(\rho_\nu(b),a)},
\where \cr
&\rho_\nu(b) = \om_b^{-1}(\rho_\nu) = \rho_\nu -
\sum_{ (b,\al)\le 0} \al, \ \
{\al\in R_+,\ \nu_\al=\nu}.
&\eqnu
\label\rhob\eqnum*
}
$$

The theorem  results immediately
in the orthogonality of $\{e_b\}$  for pairwise
distinct $b$. Macdonald also gives the formula for the squares
of $e_b$ (for $t_\nu=q^k,\ k\in \Z_+$)
and writes that he deduced it from the corresponding
formula in the  $W$-symmetric case ([C2]).
 A direct simple proof (based on the duality)
will be a subject of the next section.
Now we come to the connection between $e_b$
for $b$  from the same $W$-orbit.

Let us fix a set $\{\vep={\vep_\nu\in \{\pm 1}\}\}$
ensuring the condition  $\vep_\nu=1$ if $\nu=\nu_j$ for at least one
index $j$ such that $s_j(b_-)=b_-$. We keep the same notations
$\vep_a, \vep_j$ as for $t$.
We introduce "the $\vep$-intertwiners" (see e.g. [C2,C4]) as follows:
$$
\eqalignno{
&\Phi_j^{(\vep)}\ =\
\Bigl(T_j + (t_j^{1/2}-t_j^{-1/2})(Y_{a_j}^{-1}-1)^{-1}
\Bigr)(\phi_j^{(\vep)})^{-1},\ \ \phi_j^{(\vep)}\ =
&\eqnu
\cr
\vep_j t_j^{\vep_j/2} +
(t_j^{1/2} &-t_j^{-1/2})(Y_{a_j}^{-1}-1)^{-1} =
\vep_j(t_j^{1/2} + (t_j^{1/2}-t_j^{-1/2})(Y_{a_j}^{-\vep_j}-1)^{-1})
\label\Phi\eqnum*
}
$$
for $1\le j\le n$. They belong to the proper localization
of the affine Hecke algebra $\h_Y$ and satisfy the same relations
as $\{s_j\}$ do. Hence
$$
\Phi_w^{(\vep)}\ =\ \Phi_{j_l}^{(\vep)}\cdots\Phi_{j_1}^{(\vep)}
\for w=s_{j_l}\cdots s_{j_1}\in W
\eqnu
\label\Phiprod\eqnum*
$$
are well defined and $\Phi$ is a homomorphism of $W$.

\proclaim {Proposition}
Let  $b\in B,\ w=\om_b:\  b\to b_-$  (Proposition \ref\PIOM).
Then (see (\ref\Yone)):
$$
e_{b_-} \ =\ \prod_{(a^\vee,b)>0}
\Bigl( \vep_a {t_a-  x_a^{\vep_a}(\#b)\over
1- x_a^{\vep_a}(\#b) }
\Bigr)
\Phi_{w}^{(\vep)}(e_b),\
a\in R^\vee_+.
\eqnu
\label\Phie\eqnum*
$$
\endproclaim
\label\PHIE\theoremnum*
{\it Proof.}
The intertwiners permute $Y$:
$\Phi_{w}^{(\vep)}Y_a\ =\ Y_{w(a)} \Phi_{w}^{(\vep)}$.
Hence the r.h.s of (\ref\Phie) is proportional
to $e_-$. We put the denominators on the right
using $G$ from (\ref\G) and  formula (\ref\Phiprod):
$$
\eqalignno{
\Phi_w^{(\vep)}(e_b)\ = \
w\ &\Bigl( G_{a^l} + s_{a^l}(t_{a^l}^{1/2}-t_{a^l}^{-1/2})
(x_{a^l}(\#b)^{-1}-1)^{-1}\Bigr)
\cr
\cdots &\Bigl(G_{a^1} + s_{a^1}(t_{a^1}^{1/2}-t_{a^1}^{-1/2})
(x_{a^1}(\#b)^{-1}-1)^{-1}\Bigr)\cr
&\prod_{a\in (\la(w))^\vee}
 \vep_a {1-x_{a}^{-\vep_a}(\#b)\over
t_a^{-1/2}- t_{a}^{1/2} x_{a}^{-\vep_a}(\#b) }e_b, \where
&\eqnu
\label\Phieprod\eqnum*
\cr
\{a^1=a_{j_1},&a^2=s_{\al_1}(a_{j_2}),
a^3=s_{\al_1}s_{\al_2}(a_{j_3}),\ldots,\}=(\la(w))^\vee.
}
$$
 Due to Proposition \ref\PIOM and formula
(\ref\gonx), the conditions $a^i>0, (b,a^i)>0$
give that each quantity $\Bigl(G_a+\ldots\Bigr)$ acts
as the corresponding $t_a^{-1/2}$ on the leading $x$-term.
So the coefficient of $x_{b_-}$ in (\ref\Phieprod)
 equals
$$
\eqalignno{
& \prod_{a\in (\la(w))^\vee}\vep_a
(1-x_a^{-\vep_a}(\#b))/
(1-t_a x_a^{-\vep_a}(\#b) ).
}
$$
\proofbox

{\it Symmetric polynomials.}
The above formula results in the following explicit
expressions for the $\vep,t$-{\it symmetrizations} of $e_{b_-}$.
To introduce them we need
$$
\eqalign{
&\p_{\vep}^t\ =\ \sum_{w\in W}
\prod_\nu(\vep_\nu t_\nu^{1/2})^{\vep_\nu l_\nu(w)} T_w,
\cr
&\f_{\vep}\ =\
\sum_{w\in W}
 \Phi_w^{(\vep)}.
}
\eqnu
\label\symmetr\eqnum*
$$
Let us
 check (see (\ref\tonx) and [C2], Proposition 4.6) that
$\f_{\vep}$
is divisible on the left
by  $\p_{\vep}^t$ (i.e. $\f=\p(\ )$) and, moreover,
 $\p_{+}^t$  is divisible
by $\p_{+}\equal \p_{+}^{t=1}$ (cf. Corrolary 4.7, ibidem).
We donote the constant set $\{\vep=\pm\}$
by $\pm$.

Using $\vph$ from Section 2 (there is also a straightforward
way via the induced representations
of $\h_Y$),
$$
\eqalignno{
\vph(\Phi_j^{(\vep)})\ &=\ s_j\for \vep_j=+1,\cr
&=\ - { t_j^{1/2}X_{a_j}- t_j^{-1/2} \over t_j^{-1/2}
X_{a_j}- t_j^{1/2} } s_j
\for \vep_j=-1
&\eqnu
\label\vphPhi\eqnum*}
$$
in the basic representation.
Hence $\vph(\Phi_j^{(\vep)})+ 1$ is divisible on the right
by $T_j+ \vep_j t_j^{-\vep_j/2}$. Applying $\vph$ one more time
we get the required.

\proclaim {Proposition}
The polynomial
$$
\eqalignno{
p_{b_-}^{(\vep)}\equal &
\sum _{b\in W(b_-)}\
 \Phi_{\om_b}^{(\vep)}(e_{b_-})\ =
&\eqnu
\label\symp\eqnum*\cr
=\ &\sum _{b\in W(b_-)}
\Bigl(\prod_{(a^\vee,b)>0, a\in R^\vee_+}
\vep_a{t_a- x_a^{\vep_a}(\#b)\over
1- x_a^{\vep_a}(\#b) } \Bigr) e_b
%\eqnu
\label\psym\eqnum*
}
$$
is $\vep,t$-symmetric, i.e.
 $\p_{\vep}^t(p)$ is proportional to $p$.
 Moreover, it is $W$-invariant when $\vep=+$.
\endproclaim
\label\PSYM\theoremnum*

{\it Proof.}
Because
$\vep_\nu=1$ for $\nu=\nu_j $ if $s_j(b_-)=b_-$
for some $j$,
$\f_{\vep}(e_{b_-})$ is proportional to
$p_{b_-}^{(\vep)}$. Then it is necessary just to use
(\ref\Phie) divided by $\Phi_w$.
\proofbox

We note that the representation of $\h_Y$ in the
Macdonald polynomials for $b\in W(b_-)$ is isomorphic
to $V(b_-)=\Si(b_-)/\Si_+(b_-)$
(see (\ref\vbmin).  The image  of
$p_{b_-}^{(\vep)}$ is $m_{b_-}$ (which provides another way
to fix it uniquely).

The $+$-symmetrizations  $p_{b_-}^{(+)}$
are the Macdonald polynomials  [M1,M2].
More precisely, he defined a basis
$\{p_b, b\in B_-\}$  in the space
$\C_{q,t}[x]^W$ of all $W$-invariant polynomials
by the  conditions
$$
\eqalignno{
&p_b-m_b\ \in\ \Si_+(b),\
\langle p_b, m_{c}\rangle = 0, \when c\succ b_+,
&\eqnu
\label\macdsym\eqnum*
}
$$
for the monomial symmetric functions $m_b=\sum_{c\in W(b)}x_c$.
 One can
also introduce $\{p\}$ as eigenvectors for the
 ($W$-invariant) operators $L_f,\ $  $f\in \C_{q,t}[x]^W$:
$$
\eqalignno{
&L_{\bar{f}}(p_b)=f(q^{b}t^{-\rho} ) p_b,\
b\in B_-.
&\eqnu
\label\Lf\eqnum*
}
$$

Applying  any elements from $\h_Y$
to   $e_c (c \in W(b_-))$
we get solutions
of (\ref\Lf), because
 symmetric $Y$-polynomials are central
in $\h_Y$ (due to I. Bernstein). Since $p_{b}^{(\vep)}$
are  of this kind, Proposition \ref\PSYM
readily gives the coincidence
 $p_{b}=p_{b}^{(+)}\ (b=b_-)$ and the  coefficients
of the decomposition of $p_b$ in terms of
 $e_c, c\in W(b)$.
The formula for these coefficients was
announced in [M3] (where
$t_\nu=q^{k}, k\in Z_+$). In the differential case,
the coefficients (for arbitrary $\vep$) were calculated
in [O2].

The elements $p_{b}^{(-)}, b\in B_-,$ are also quite
remarkable. The map
$$
\eqalign{
&p_{b}\to p'_b\ =\ p_{b}^{(-)}/\hbox{det}_t \for\cr
&\hbox{det}_t\ =\ \prod _{a\in R^\vee_+}((t_a x_a)^{1/2}-
(t_a x_a)^{-1/2})
}
$$
is exactly the action of the shift operator from [C2].
Namely, $p'_b$ is proportional to the
Macdonald symmetric polynomial for $t'_\nu=t_\nu q_\nu$,
$b'=b+\rho$.
In the differential
case this observation is due to Opdam (ibid.). It is closely
connected with the main theorem from [FV]. Macdonald also
uses this approach to the shift operators
in [M3] (the difference case).

%
%
%		Section 5
%
%
%\vskip 10pt
\section { Duality, applications}
First of all we will use Theorem \ref\dual to define the
{\it
Fourier pairing}. In the classical theory the latter is the
 inner product of a function and the
Fourier transform  of another function.
In this and the next sections we will
continue to identify the elements
$H\in \HH\ $ with their images $\hH$. The following pairing
on $f,g\in \C_{q,t}[x]$ is symmetric and non-degenerate:
$$
\eqalign{
 &[\![f,g]\!]= [\![f(X),g(X)]\!] = [\![\vph(f(X))g(X)]\!]  =\cr
&[\![\bar{f}(Y)g(X)]\!] =  \{L_{\bar{f}}(g(x))\}(t^{-\rho}),\cr
&\bar{x}_b\ =\ x_{-b}\ =\ x_b^{-1},\ \bar{ q}\ =\  q,\
\bar{t}\ =\ t,
}
\eqnu
\label\Fourier\eqnum*
$$
where
$L_f$ is from Theorem \ref\YONE ,
and we used the main defining property (\ref\hat)
of the representation
from  Theorem \ref\faith.

The Fourier adjoint $\vph(L)$ of any $\C_{ q,t}$-linear
operator $L$
acting  in $\C_{ q,t}[x]$ is defined from the relations:
$$
\eqalign{
&[\![L(f),g]\!]\  =\   [\![f,\vph(L)(g)]\!],\ f,g\in \C_{ q,t}[x].
}
\eqnu
$$
This anti-involution  ($\vph^2=\hbox{id}$) extends  $\vph$ from
(\ref\vph) by construction. If $f\in C_{ q,t}[x]$, then
$\vph(L_f)= \bar{f}(X)$. We arrive at the following
theorem:

\proclaim {Duality Theorem}
Given $b,c\in B$ and the corresponding Macdonald's
polynomials $e_b, e_c$,
$$
\eqalign{
&e_b(\#c)e_c(\#)\ =\ e_b(t^{-\rho(c)} q^{c})e_c(t^{-\rho})\ =\
[\![e_b,e_c]\!]\ =\cr
& e_c(\#b)e_b(\#)\ =\ e_c(t^{-\rho(b)}q^{b})e_b(t^{-\rho})\ =\
[\![e_c,e_b]\!]
}
\eqnu
\label\pp\eqnum*
$$
in the notations from (\ref\rhob).
\endproclaim
\label\DUAL\theoremnum*
\proofbox

To complete the theorem we will calculate $e_b(t^{-\rho})$
 together with the norms $\langle \ep_b,\ep_b\rangle$ of
the {\it renormalized Macdonald polynomials}
$\ep_b\equal e_b/e_b(t^{-\rho})$ by means of the Recurrence Theorem.

\proclaim{Main Theorem}
Given $b\in B$, let   $b_-,\ b_+=w_0(b_-)$
be the corresponding elements from $W(b)\cap B_{\pm}$,
$b^o=-w_0(b).$
Then
$$
\eqalignno{
&e_b(t^{-\rho})\ =\ x_{b_-}(t^{\rho})\prod_{a\in R_+^\vee}
\Bigl(
{
1- q_a^{j}t_a x_a(t^\rho)
 \over
1- q_a^{j}x_a(t^\rho)
}
\Bigr),
&\eqnu
\label\EV\eqnum*\cr
&\langle \ep_b,\ep_b\rangle\ =\
\prod_{a\in R_+^\vee}
\Bigl(
{
t_a^{1/2}-q_a^jt_a^{-1/2} x_a(t^\rho)\over
t_a^{-1/2}-q_a^jt_a^{1/2} x_a(t^\rho)
}
\Bigr),
&\eqnu
\label\NORMS\eqnum*
}
$$
where the products are over the same set
$J_b = \{j \}$:
$$
\eqalignno{
& 0<j<(a^\vee,b_+)\if (a,b^o)>0,\ \
 0<j\le(a^\vee,b_+)\if (a,b^o)<0.
&\eqnu
\label\jbset\eqnum*
}
$$
\label\MAIN\theoremnum*
\endproclaim
\proofbox

 We mention that there is a straightforward passage
 to non-reduced root systems and
to  $\mu$  introduced for
$\al\in R_+$ instead of $a\in R_+^\vee$ (see [C2]).
As to the latter
case, it is necessary just to replace the indices $a$ by $\al$
($ q_a \to q, \rho\to r$) in the formulas for $\{e_a\}$.
In the $W$-symmetric case these statements
(the Macdonald conjectures) are from [C2,C3].
In [M3], Macdonald gives a formula for the norms of $e_b$.
Hopefully it coincides with (\ref\NORMS) after the multiplication
by $e_b(t^{-\rho})(e_b(t^{-\rho}))^*$ and then by $\langle \mu\rangle$.
 In his paper,
$t_\nu=q^k,k\in\Z$. The differential case is due to Opdam.

{\it Discretization.} Let us establish the {\it
recurrence relations} for the Macdonald polynomials
generalizing the three-term relation for the $q$-ultraspherical
polynomials (Askey, Ismail) and the Pieri rules.
We follow [C3] where the symmetric case was considered.
We need to go to the lattice
version of the  functions and operators. The
 discretization of
functions $g(x)$ in $x\in \C^n$  and the operators
acting on such functions
is defined
as in Theorem \ref\YONE :
$$
\eqalignno{
&{}^\de x_a(bw)=x_a(q^b t^{-w(\rho)})=
q^{(a,b)}\prod_\nu t_\nu^{-(w(\rho_\nu),a)},\cr
& ({}^\de \hu ({}^\de g))(bw)\ =\ {}^\de g(\hu^{-1}bw),
\ \hu\in W^b,
&\eqnu
\label\deltaf\eqnum*
\cr
&({}^\de X_a({}^\de g))(bw)\ =\  x_{a}( bw)
\ {}^\de g(bw).
}
$$
It is a homomorphism. The image is the space
$\hbox{Funct}(W^b, \C_{q,t})$ of
 functions on $W^b$ and   operators acting on
such functions. We will sometimes omit
${}^\de$ and put $g(\hw)$ instead of  ${}^\de g(\hw)$
etc.

Given an arbitrary linear combination of functions
$\{\phi_{\hw}(\ ), \hw\in W^b\}$, we can also apply the above operators
to the sufficies:
$$
\eqalign{
& {}_\de (g(x)\hu)(\sum_{{\hw}\in W^b}c_{\hw}\phi_{\hw}(\ ))\ =
\ \sum_{{\hw}\in W^b}c_{\hw} g( \hw)\phi_{\hu^{-1}\hw}(\ ),\
c_{\hw}\in \C.
}
\eqnu
\label\deltsuf\eqnum*
$$
It is an anti-homomorphism,
 i.e.
$$
{}_\de (GH)\ =\  {}_\de H\ {}_\de G \for \hbox{\
operators\ } G,H.
$$

We will mostly use the discretizations
$ \ep_b({\hw})= e_b ({\hw})/e_b (0)$
of the
renormalized Macdonald polynomials
 $\ep_b(x)= e_b(x)/e_b(t^{-\rho})$, and
especially
$\ep_b(\#c)=\ep_c(\#b),$ where
$\#c\equal\pi_c=c\om_c^{-1}.$
See (\ref\Yone),(\ref\rhob),
and Theorem \ref\DUAL.
Sometimes we drop $\#$ and write $\ep_b(c)$ instead of
$\ep_b(\#c)$. For example, $\ep(0)$ always means $\ep(\#)$.
Vice versa, we will consider  the sufficies $b$ as
 elements from $W^b$
via the same  map $b\to \#b$.

Given a  polynomial $f\in \C[x]$, we
construct the operator $L_f=f(Y)$, go to its discretization
${}^\de L_f$,
and finally introduce the {\it recurrence operator}
 $\La_f= {}_\de L_f$  acting on the sufficies
${\hw}\in B$ of any $\C$-valued
functions
$\phi_{\hw} (\ )$. We write $\La_a$ when $f=x_a, a\in B.$

\proclaim{Recurrence Theorem}
For arbitrary $a,b\in B, f\in \C[x]$,
$$
\eqalign{
&\La_{f}(\ep_b(x))\ =\  \bar{f}(x)\ep_b(x),\
\ \La_{a}(\ep_b(x)) =  x_a^{-1} \ep_b(x),
}
\eqnu
\label\pbpb\eqnum*
$$
where  $\bar{f}(x)=f(x^{-1})$. The operators $\La$ (acting
on $\#b$) are well defined. It  means that they do
 not produce  the indices which do not belong to
 $\#B=\{\pi_b,b\in B\}$.
\label\RECUR\theoremnum*
\endproclaim

{\it Proof.}
We can rewrite (\ref\Lf) as follows:
$$
\eqalignno{
&{}^\de L_f({}^\de e_b)\ =\ \bar{f}(\#b)\ {}^\de e_b.
&\eqnu
\label\deLf\eqnum*}
$$
Replacing $e$ by $\ep$ and using the duality we yield:
$$
\eqalignno{
{}^\de & L_f(\ep_b(c))\ =\ \bar{f}(\#)
\ {}^\de \ep_b(c),\cr
& \La_f(\ep_b(c))\ =\ {}^\de\bar{f}(c)
\ {}^\de \ep_b(c),&\eqnu
\label\deLaf\eqnum*
}
$$
and (\ref\pbpb) if we can ensure that  $\La_f$ does not
create polynomials $\ep_c$ with the indices apart from $\#B$. The
latter will
be checked  in the next section.
\proofbox

The theorem has many applications. For instance, we can
prove the Main Theorem.
To demonstrate this let (see Proposition \ref\YBGF)
$$
\eqalign{
& Y_a  = \sum_{b\succeq a, w\in W}  (bw)  g_{a}^{bw},\
\La_a  = \sum_{b\succeq a, w\in W} {}_\de (g_{a}^{bw}){}_\de (bw),
\cr
 & Y_{-a}  = \sum_{b\succeq a, w\in W}   f_{a}^{bw}(bw)^{-1} ,\
\La_{-a}  = \sum_{b\succeq a, w\in W}{}_\de (bw)^{-1} {}_\de (f_{a}^{bw}).
}
\eqnu
\label\ella\eqnum*
$$
Thanks to the theorem:
$$
\eqalign{
& \La_a \ep_b \ = \
x_a^{-1}\ep_b\ =\
\sum_{\hw}  g_{a}^{\hw}(\hw^{-1}\#b)\ep_{\hw^{-1}\#b},
\cr
& \La_{-a} \ep_b \ = \
x_a\ep_b\ =\
\sum_{\hw}  f_{a}^{\hw}(\#b)\ep_{\hw\#b}.
}
\eqnu
\label\Laep\eqnum*
$$
For example,
$$
g_a^{0}(\#)\ =\ \langle\mu_1 x_a^{-1}\rangle\ =\
\langle\mu_1 x_a\rangle^* \ =\ f_a^{0}(\#)^*
$$
are the coefficients of $\mu_1$ (we remind that
$g(\#)=g(t^{-\rho})$). Their   description is one of the main
open problems in the Macdonald theory (we hope to consider it in
the next papers).

 \vfil
\proclaim{Proposition} The presentation $x_a\ =\ \sum_{b\preceq a}
f_a^{\#b}(\#)\ep_{b}$ leads to:
$$
\eqalignno{
&f_{a}^{\#b} (\#)^*\langle \ep_b,\ep_b\rangle\ =\
\langle \ep_b,x_a\rangle\ =\
g_a^{\#b}(\#),
&\eqnu
\cr
\label\fqgq\eqnum*
&e_b(\#)\ =\ f_b^{\#b}(\#),\
\langle \ep_b,\ep_b\rangle\ =\
g_b^{\#b}(\#)(f_b^{\#b}(\#)^*)^{-1}.
&\eqnu
\label\norms\eqnum*
}
$$
\label\FMG\theoremnum*
\endproclaim

{\it Proof.}
The first relation results from (\ref\Laep) for $b=0$.
It gives readily the formula for $e_a(\#)$. Then
$$\langle \ep_b,x_a\rangle = \langle \ep_b,
\sum_{b} f_a^{\#b}(\#)\ep_{b}\rangle =
(f_a^{\#b}(\#))^*\langle \ep_b,e_b\rangle.
$$
On the other hand,
$$\langle \ep_b,x_a\rangle\ =\
\langle \mu_1\ep_b x_a^{-1}\rangle\ =\ g_a^{\#b}(\#),
$$
which is (\ref\fqgq).
Letting $a=b$, we come to the last formula for the norm.
\proofbox

{\it Proof of the Main Theorem.}
In fact, the coefficients $f_b^{\#b}(\#),\ g_b^{\#b}(\#)$
were  calculated in Proposition \ref\YBGF.
We need only to substitute the evaluation of $T_{\om_b}$
 at $\#$, that is $l_t(\om_b)=\prod_\nu t_\nu^{l_\nu(w)/2}$.
\proofbox

%
%
%		Section 6
%
%
%\vskip 10pt
\section {Roots of unity}
Let us assume that $ q$ is a primitive
$N$-th root of unity for $N\in \N$
and first consider $t$ as an indeterminate parameter.
More precisely, we will operate over the field
$\Q_t^0\equal \Q( q_0,t)$ where we fix $ q_0$
such that $ q_0^{2m}= q$ ($ q_0$ belongs to a proper
extension of $\Q$). Actually all formulas will  hold
even over the localization of $\Z[ q_0,t]$
by $t^{r_1} q^{s_1}(1-t^{r_2} q^{s_2})\neq 0$, $r_i,s_i\in \Z$.

The pairing
$$
\eqalign{
B\times B\ni a\times b\to  q^{(a,b)}\equal  q_0^{2m(a,b)}
}
\eqnu
\label\BtimeB\eqnum*
$$
acts through $B_N\times B_N$, where
 $B_N\equal B/K_N$,  $K_N$ is its radical.

Following the previous section
we restrict the functions $\{x_b\}$ and
the operators $\{Y_b\}$ to the
 $W^b$ using the pairing (\ref\BtimeB) and
the formulas $x_a(bw)=x_a(q^bt^{-w(\rho}),\ b\in B,w\in W$.
Given $w\neq u$ and any $b,c$, there exists $a\in B$
such that
 $x_a(bw)\neq x_a(cu)$ ($t$ is generic). Hence the discretization
maps via $W^b_N\equal B_N\rsmash W$.

The $T,Y$-operators
are well defined over $\Q_t^0$ since their
 denominators
are  products of the binomials $(x_a q^k-1)$
for $a\in R^\vee, k\in \Z$. The latter remain
non-zero when evaluated at $ q^b t^{-w(\rho)}$
since $(a,\rho)$ never equals $0$ ($x_a(t^{-w(\rho)}$
always contains $t$). Hence the discretizations of these
operators exist too.
More exact information
about the properties of these coefficients can
be extracted from Proposition \ref\FMG.

Let $B(N)\subset B$ be a fundamental domain of the
group $K_N$. It means that the map
$B(N)\ni b\to \#b\in W^b_N$
is an isomorphism.
Further we identify these two sets, putting
$$
B(N)\ =\
\{\be^1,\ldots, \be^d \}\  =\ W^b_N , \where d=|W^b_N|,\ \be^1=0.
\eqnu
\label\bei\eqnum*
$$
The images of  $\hw\in W^b$
in $B(N)$ will be denoted by $\hw'$
(i.e. $\hw'=\be^i$ if $\hw=\#\be^i\hbox{\ mod\ }K_N$.
One may assume that
$-w_0(B(N))=B(N)$ for the longest element $w_0$.

Abusing the notations,  we write $g(b)$ where $b\in B$ instead of
${}^\de g(\#b)$. Correspondingly,
by $b'$ we mean the image of $\#b$ in $B(N)$, that is
the image of $b$ in $B_N$.

Let us consider (temporarily) the case when
$N$ is coprime with the order $ |B/A|= |O|$ taking
$ q_0^2= q^{1/m}$ in the $N$-th roots of unity.
Then $K_N=  NP\cap B$ for the weight lattice
$P=\oplus_{i=1}^n\Z\om_i$ generated by the  $\om_i$
(dual to $a_i$). We can take the following fundamental domain:
$$
\eqalign{
 &B(N)\subset  \{c\in W(b_-),\where
b_-= -\sum_{i=1}^n k_ib_i\in B_-\}, \cr
&0\le k_i\le N \if (2/\nu_i, N)=1,\ \ 0\le k_i\le\nu_i N/2
\hbox{\ \ otherwise},
}
\eqnu
\label\kineq\eqnum*
$$
removing $b_-$ such that $ k_i=0,k_j= N/(2/\nu_i, N)$
for at least one pair of the indices. Moreover, if
$c_-\in a_-+NP_+$ then we do not take $a\in W(a_-)$ when
$\om_a\in \{\om_c,c\in W(c_-)\}$.
Recall  that $2/\nu_i=2/(\al_i,\al_i)= 1,2,3,\ \#b=b\om_b^{-1}$
(see Section 1).

Let us demonstrate that the Macdonald polynomials $e_b$
are well defined for $b\in B(N)$ (later we will
see that they  always exist). We  introduce them
 directly from (\ref\Lf), using that the $Y$-operators
preserve any subspaces
$$
\Si_*^0(b) = \oplus_{ c\succ b} \Q_t^0 x_c,\and
\Si^0(b)= \Si_*^0(b)\oplus \Q_t^0 x_b,\ b,c\in B.
$$

It is necessary to check that given $B\ni c\succ b$,
there exists at least one $ a\in B$ such that
$x_a(q^{b}t^{-\rho(b)})\neq$ $  x_a(q^{c}t^{-\rho(c)})$
for
$$\rho(b) \ = \ \om_{b}^{-1}(\rho)\ =\
\rho -
\sum \al,\ \ {\al\in R_+, \ (b,\al)\le 0}
$$
(see (\ref\rhob)).
Then the eigenvalues of $Y$ will separate $e_b$ from
the elements from $\Si_*^0(b)$ and we can argue by induction.

 It is obvious if $\om_b\neq\om_c$.
Otherwise
$ k_i< N (2/\nu_i, N)^{-1}$ and $c_-\succ b_-$.
Then we repeat the corresponding reasoning from
[C3], Section 5.

 The discretizations $\ep_{b'}(\hw)$
are well defined too and
depend only on the images $\hw'$ because
$\ep_{b'}$ are linear combinations of $x_a, a\in B$.
We see that
 $\{\ep_{\be^i}(\hw')\}$ form a basis in the
space $V_N\equal \hbox{Funct}(B(N), \Q_t^0)$
of all $ \Q_t^0$-valued
 functions on $B(N)$. Indeed, they are non-zero and
the action of the $\{{}^\de Y_a\}$ ensures that they are linearly
independent.

{\it The end of the  proof of Theorem \ref\RECUR }.
First of all, let us rewrite formally relation
 (\ref\pbpb) for $f=x_a$
 as follows:
$$
\eqalign{
&x_a^{-1}\ep_b(x)\ =\ \La_a^{\#} (\ep_b(x))+
 \sum_{\hw\not\in \#B} M_{ab}^{\hw} \ep_{\hw}(x).
}
\eqnu
\label\mpie\eqnum*
$$
Here $\hw$ form a finite set $E=E(a,b)\ (E\cap \#B=\emptyset)$,
$M_{ab}^{\hw}$ are rational functions of $ q,t$.
The
truncation $\La_a^{\#}$
of $\La_a$ is uniquelly determined by the
 condition that it does not contain $\hw$ moving $\#b$ to
elements apart from $\#B$. Assuming that $N$ is sufficiently
big the discretization gives the relation (see (\ref\deLaf)):
$$
\eqalign{
&x_a^{-1}(c)\ep_b(c)\ =\ \La_a^{\#} (\ep_b(c))+
 \sum_{\hw\not\in \#B} M_{ab}^{\hw} \ep'_{\hw}(c),\  c\in B(N),
}
\eqnu
\label\mpiede\eqnum*
$$
for $ x_a(c)={}^\de x_a (c)$.
Here   $\ep'_{\hw}=\ep_{\hw'}$ for $\hw\in E$. This substitution was
impossible before the discretization. We remind
that  the  formula  with
 $\ep_{c}(\hw)$ in place of  $\ep_{\hw}(c)$ is always true.
Because $c$ is taken
from $B(N)$ the  discretization of $\ep_c$
 exists. Therefore  we can replace the argument $\hw$ by
 $\hw'\in B(N)\subset B$,
 and then go from  $\ep_c(\hw')$ to $\ep'_{\hw}(c)$
thanks to the duality.

As to $M_{ab}^{\hw}$, they  are  the values of the coefficients
of ${}^\de Y_a$ and are  well defined when
$ q^N=1$  (enlarging $N$ we can get rid of singularities
in $ q$ even if $M$ are arbitrary rational).

On the other hand :
$$
\eqalign{
&x_a^{-1}\ep_b(x)\ =\
 \sum_{h\in B} K_{ab}^{h} \ep_{h}(x),
}
\eqnu
\label\mpix\eqnum*
$$
where the coefficients $K_{ab}^h$ are rational functions
of $ q,t$,
$\{h\}$ form a finite set $H =H(a,b)\subset B$. The discretization
 gives that
$$
\eqalign{
&x_a^{-1}(c)\ep_b(c)\ =\
 \sum_{h\in B} K_{ab}^{h} \ep_{h}(c), \ c\in B(N).
}
\eqnu
\label\mpixde\eqnum*
$$
We pick  $N$ to avoid possible singularities.

Since $N$ is sufficiently big, the eigenvalues of the
$Y$-operators distingwish all $\ep_{d}(c)$  for $d\in
(E)'\cup H$. It holds only for generic $t$ (say, when
$t=1$ it is wrong). Comparing (\ref\mpiede) and
(\ref\mpixde) we conclude that $M_{ab}^{\hw}=0$ for all $\hw\in E$,
when $ q^N=1$.
Using again that $N$ is arbitrary (big enough, coprime with $|O|$)
we get
that the actions
of $\La_a^{\#}$ and $\La_a$   coincide on $\ep_b$, i.e.
the latter operator does not create the indices  not
from $\#B$.
\proofbox

Let us go back to the general case (we drop the condition
$(N,|O|)=1$).
Once the Recurrence Theorem has been established we can
use Proposition \ref\FMG without any reservation.
It readily gives that the Macdonald polynomials $e_b,\ep_b,$
and  ${}^\de\ep_b$ are well defined for
arbitrary $b\in B$ because
 $f_b^{\#b}(t^{-\rho})\neq 0$ (formulas (\ref\EV) and (\ref\norms)).
 Moreover,
${}^\de\ep_b = {}^\de\ep_c$ if and  only  if
$b'=c'$. Hence
the {\it restricted Macdonald polynomials} $\ep_{\be^i}(\tw),
1\le i\le d$ (see (\ref\bei)) form a basis in
 $V_N= \hbox{Funct}(B(N), \Q_t^0)$.
Indeed,   $\ep_{\be^i}(c')$ are eigenvectors of the
 ${}^\de L$-operators separated by the eigenvalues. They
 are always non-zero since
$\ep_b(0)=1$. Hence they are linearly independent over
$\Q_t^0$ and form a basis in $V_N$. Every ${}^\de\ep_b$
is an $L$-eigenvector and coincides
with  one of them (when $\be^i=b'$).

 Similarly,
${}^\de L_{\ep_b} = {}^\de L_{\ep_c}$ if and  only  if
$b'=c'$ because the latter condition is necessary and
sufficient to ensure the coincidence of the sets of
eigenvalues.
We will also use the basis of
the {\it delta-functions}
$\de_{\be^i} (\be^j)\equal\de_{ij}$ separated
by the action of $\{{}^\de x_a\}$.

\proclaim {Proposition}
 The discretization map
 supplies
$V_N$ with the structure of
an $ \HH$-module which is irreducible.
 The Fourier pairing is well defined
on $V_N$ and induces the anti-involution
$\vph$.
\endproclaim
\label\IRRED\theoremnum*
{\it Proof.}
 The radical of the Fourier pairing (\ref\Fourier)
 contains the kernel of the discretization
map $\Q_t^0[x]\to V_N$. Its restriction
 to $V_N$ is  non-degenerate since
 $$
\eqalign{
\Pi=(\ep_{ij}), \where \ep_{ij}\ =\
[\![\ep_{\be^i},
 \ep_{\be^j}]\!]= \ep_{\be^i}(\be^j),
}
\eqnu
\label\pimat\eqnum*
$$
is the matrix connecting
the bases $\{\ep\}$ and $\{\de\}$.
The coresponding anti-involution coincides with
 $\vph$.
Thus  $V_N$
 is semi-simple.

If $V_N$ is reducible then
$\ep_{\be^1}=1$   generates a
proper $\HH$-submodule ($\neq V_N$). But it  takes
 non-zero values at any points of $B(N)$. Hence
its $\HH$-span must contain all $\de_{\be^i}$. We come to
a contradiction.
\proofbox

{\it When $t$ are roots of unity.} Till the end
of the paper  $t_\nu= q_\nu^{k_\nu}$ for
$k_\nu\in \Z_+,\  \nu\in \nu_R$. The $L$-operators
act in  $\Q^0[x]$ for $\Q^0\equal \Q( q_0)$.
Let $J\subset \Q^0[x]$ be the radical of the pairing
$[\![\ ,\ ]\!]$. It is an ideal and an
$\HH$-submodule. The quotient (a ring and an $\HH$-module)
$\v=\Q^0[x]/J$ is finite dimensional over $\Q^0$.
It results from Proposition \ref\IRRED as generic
$t$ approaches $ q^k$.

The set $\tB$ of maximal ideals
of $\v$ will be considered as a subset of $B_N=B/K_N$
relative to the map:
$$
b\ \mapsto\ q^{\#b},\  \#b\ =\ b-\om_b^{-1}(\rho_k),\
\rho_k\equal \sum_\nu k_\nu\rho_\nu.
\eqnu
\label\rhok\eqnum*
$$
It  contains $0$  corresponding to the evaluation
at $q^{\#}$, since
$[\![1,f]\!]=f(0)$ and  $J$ belongs to the ideal
$\{f,f(0)=0\}$.
We put
$$
\tB\ =\
\{\tbe^1,\ldots, \tbe^\pa \}\  \subset\  B, \ \tbe^1=0.
\eqnu
\label\tbei\eqnum*
$$
We keep the same abbreviations:
$$
\eqalign{
&^\de g(bw)\ =\ g( q^{b-w(\rho_k)}),
\hbox{\ \ however\ \ } g(b)\ =\ g(\#b).
}
\eqnu
\label\gbwgb\eqnum*
$$

By the construction,   all $\{Y_a\}$-eigenvectors
in $\v$ have pairwise distinct eigenvalues (the difference
of any two of them  with the same sets of
eigenvalues and coinciding evaluations at $0$
belongs to $J$). Applying this to  $\ep_0=1$ generating
$\v$  as an $\HH $-module
we establish the irreducibility of $\v$ (use the pairing
$[\![\ ,\ ]\!]$ and follow Proposition \ref\IRRED).

Next, we will introduce the {\it restricted Macdonald
pairing} (cf. (\ref\mu),(\ref\innerpro)):
$$
\eqalign{
&\langle f(x),g(x)\rangle'\equal\sum_{c\in B_N}
\mu(c)f(c)\bar{g}(c) \for f,g\in \Q^0[x]^W,\cr
&\mu\ =\ \prod_{a\in R^\vee_+}
(1-x_a q_a^{k_a-1})\cdots(1-x_a)\cdots
(1-x_a^{-1} q_a^{k_a-1})(1-x_a^{-1} q_a^{k_a}).
}
\eqnu
\label\resinner\eqnum*
$$
Here $\mu(c)={}^\de\mu(c)=\mu(\#c)$.

The same verification as in [C2],  Proposition 4.2
gives that
$$
\langle Af,g\rangle ' \ =\ \langle f,A^*g\rangle'
\eqnu
\label\resinv\eqnum*
$$
for the anti-involution ${}^*$ from (\ref\star)
considered on  the operators
 from $\HH$ acting on  polynomials.

Let us assume that:
$$
 q_a^{(\rho_k,a)+i}\neq 1 \hbox{\ for\ all\ } a\in R^\vee_+.\
i=-k_a+1,\ldots,k_a,
\eqnu
\label\krho\eqnum*
$$
 In the simply
laced case ($A,D,E$), it is equivalent to the condition
 $N> k((\rho,\th) + 1)$.

\proclaim {Lemma}
The natural map $\v\to \tV \equal \hbox{Funct}(\tB,\Q^0)$
is an isomorphism
which supplies $\tV$ with the structure of a non-zero
irreducible $\HH$-module. Both
 pairings $[\![\ ,\ ]\!], \langle\ ,\ \rangle$ are
well defined and non-degenerate on $\tV$.
\endproclaim
\label\COIN\theoremnum*

{\it Proof.}
The radical $J'$ of the pairing
$\langle\ ,\ \rangle'$ in $\Q^0[x]^W$ is an $\HH$-submodule.
It equals the space of all functions $f(x)$ such that
${}^\de f(c)=0$ for $c$ from the subset $B'\subset B_N$
where $^\de\mu$ is non-zero. The  set $B'$ contains
$0$, since
$\mu'(t^{-\rho})\neq 0$ because of  condition (\ref\krho).
Hence the linear span $J+J'$ (that is an $\HH$-submodule)
does not coincide with the entire $\Q^0[x]$. and
the  irreducibility of $\v$ results in $J+J'=J$.
\proofbox

Introducing now the delta-functions $\tde_i=\de_{\tbe^i}$,
we can define the $\ep$-functions $\{\tep_i\}$
from the orthogonality and evaluation conditions
$$
[\![\tep_i,\tde_j]\!]= C_i\de_{ij} \and \tep_i(0)=1,\
, 1\le i,j\le \pa.
\eqnu
\label\tpitde\eqnum*
$$
They are eigenvectors of the $Y_a$-operators with
the eigenvalues $x_a^{-1}(\tbe^i)$ and linearly
 generate $\tV$. The sets of eigenvalues are pairwise
distinct and $\langle \tpi_i,\tpi_i\rangle'\neq 0$.

Presumably the $\ep$-functions are the discretizations of
certain restricted Macdonald's polynomials and the
above scalar products can be calculated explicitly
but we will not discuss this here.

We will  use that
$(K_N,r_\nu)\in N\Z, $ where $r_\nu=(2/\nu)\rho_\nu\in B$.
Let us impose one more restriction:
 $$
q^{(a,a)/2}=q_0^{m(a,a)}\ =\ 1   \for a\in K_N, \nu\in\nu_R.
\eqnu
\label\aarho\eqnum*
$$
If $ q_0$ is a primitive root of degree $2mN$ then
$K_N=NQ\cap B$ for the root lattice $Q=\oplus_{i=1}^n \Z\al_i$
(see (\ref\BtimeB)). This condition  obviously holds true
for even $N$ (all roots systems). For odd $N$, it is necessary
to exclude  $B_n, C_{4l+2}$. In the latter case, $B\subset Q, m=1$
and we can
pick $q_0$ in the roots of unity of degree $N$.

\proclaim {Theorem}  Introducing $\tPi=(\tep_{i}(\tbe^j))$
(see (\ref\pimat,\ref\tbei)) for $\tep_i=\ep_{\tbe^i}$,
 let
$$
\eqalign{
&\t_+ = \hbox{Diag}( q^{(\tbe^i,\tbe^i)/2}x_{\tbe^i}(t^{-\rho})),\
\t_- = \Pi\t_+^{-1}\Pi^{-1},\ \Om  =
\t_+^{-1}\t_-\t_+^{-1}.
}
\eqnu
\label\DETOm\eqnum*
$$
The conjugations be these elements induce the automorphisms
$\tau_{\pm},\om$ of \HH\ acting in $\tV$.
Let us decompose $\tV=\oplus\tV_\chi$ relative to the central characters
$\chi$  of $\H=\Q^0[T_1,\ldots,T_n]$. The following map gives
 a  projective action of the group $SL_2(\Z)$ in every $\tV_\chi$:
$$
\Bigl(\matrix 1&1\\  0&1\endmatrix\Bigr) \to \t_+,\
\Bigl(\matrix 1&0\\  1&1\endmatrix\Bigr) \to \t_-,\
\Bigl(\matrix 0&-1\\  1&0\endmatrix\Bigr)\to \Om.
$$

\label\ETD\theoremnum*
\endproclaim

{\it Proof.}
Setting
$x_b= q^{z_b}, \ z_{a+b}=z_a+z_b, \ z_i=z_{b_i},\
a(z_b)= z_b-(a,b), \ a,b\in \R^n,$
we introduce formally the {\it Gaussian}
$\ga\ =\  q^{\Sigma_{i=1}^n z_i z_{\al_i}/2}$,
which  satisfies the following (defining) difference relations:
$$
\eqalign{
&b_j(\ga)\ =\  q^{(1/2)\Sigma_{i=1}^n (z_i-(b_j,b_i))
(z_{\al_i}- q_i^j)}\ =\cr
& \ga  q^{-z_j+ (b_j,b_j)/2 }\ =\  x_j^{-1}\ga  q^{(b_j,b_j)/2}
\for 1\le j\le n.
}
\eqnu
\label\gauss\eqnum*
$$

The Gaussian  commutes with $T_j \for 1\le j\le n$
because it is $W$-invariant.
When $b_r$ are minuscule ($r\in O^*$), we  use directly  formulas
(\ref\Yb, \ref\Demaz) to check that
$$
\ga(X)Y_r\ga(X)^{-1}\ =\  X_r q^{-(b_r,b_r)/2} Y_r \ =
\ \tau_+(Y_r).
$$
A straightforward calculation gives that
$$
\eqalign{
&\ga(X)T_0\ga(X)^{-1}\ =\  \tau_+(T_0) = X_0^{-1}T_0^{-1},\cr
&\tau_-(T_0)\ =\  T_0,\
\om(T_0)\ =\ X_\th^{-1}Y_\th^{-1}T_0.
}
\eqnu
\label\gato\eqnum*
$$
Hence the conjugation by $\ga$ induces $\tau_+$.

The formula for $\t_+$ describes multiplication
by $\ga$
in $\tV$
(up to a constant factor) in the basis of delta-functions.
 Really,
$$
\eqalign{
&\ga(c)\ =\ {}^\de \ga (c)\ =\
  q^{\Sigma_{i=1}^n \ze_i \ze_{\al_i}/2}\for \cr
&\ze_i\  =\ \log_q( x_i( q^c t^{-\rho}))\ =\
(b_i,c)+\log_{q} (\prod_\nu t_\nu^{-(\rho_\nu,b_i)}).
}
\eqnu
\label\taugau\eqnum*
$$
Hence $ \ga(c) = g  q^{(c,c)/2}x_c(t^{-\rho})$ for
$g =  q^{(\rho_k,\rho_k)}.$
Since  the  matrix $\t_+$  is important up to proportionality
one can  drop the constant $g$.
We see that changing $c$ by any elements from $K_N$
does not influence $\ga(c)$ because of the condition
(\ref\aarho), which makes the multiplication by $\ga$
well defined.

Next, the automorphism $\tau_- = \vph \tau_+\vph$
corresponds to $\tPi \t_+^{-1}\tPi^{-1}$, and the matrix $\Om$
from (\ref\DETOm) induces
 $\om=$ $\tau_+^{-1}\vph\tau_+\vph\tau_+^{-1}$
 in the same delta-basis.
Indeed, $\tau_-$ is the application
of $\ga(Y)$. It multiplies $\pi_b$ by $\ga(t^\rho q^{-b})$
whereas $\ga(X)$ multiplies $\de_b$ by $\ga(t^{-\rho} q^{b})$
(so we need to  inverse  $\t_+$). More formally,
one can use the equation
$[\![\t_+^{-1}f,g]\!]=[\![f,\t_-g]\!]$.

Finally, any relations from $SL_2(\Z)$ hold for these matrices up
to proper central elements of $\H$ (Theorem \ref\GLZ). Thus the last
statement results directly from
the irreducibility
of $\tV$.
\proofbox

The theorem is a non-symmetric version of the last theorem
from [C3]. The latter in its turn
generalizes the construction due to Kirillov [Ki]
(in the case of $A_n$)
and is directly related to
Theorem 13.8 from [K]  when $t= q$.
Following [C3] one can extend the above map to the projective
$GL_2(\Z)$.  The biggest projective representations
of $SL_2(\Z)$ or $GL_2(\Z)$ can be obtained
from the eigenvalues of the element $T_{w_0}^2$ in $\tV$.
These groups act projectively
 in the corresponding spaces of eigenvectors.

Presumably the results from the last section have counterparts
for generic
$ q,t$ in the analitic setting.
They  are connected with the monodromy representation of the
double affine Knizhnik-Zamolodchikov equation and the main theorem
from [KL2]. They also might help
 to  renew
elliptic functions towards the Ramanujan theories.

%
%
%
%      REFERENCES
%
%
%
%\vskip 15pt
\AuthorRefNames [BGG]
\references
%\medskip
%\ninerm
%\baselineskip=11pt %!
\vfil

[B]
\name{N. Bourbaki},
{ Groupes et alg\`ebres de Lie}, Ch. {\bf 4--6},
Hermann, Paris (1969).

[C1]
\name{I. Cherednik},
{ Double affine Hecke algebras,
Knizhnik- Za\-mo\-lod\-chi\-kov equa\-tions, and Mac\-do\-nald's
ope\-ra\-tors},
IMRN (Duke M.J.) {  9} (1992), 171--180.

[C2]
\bibline,{ Double affine Hecke algebras and  Macdonald's
conjectures},
Annals of Mathematics {141} (1995), 191-216.

[C3]
\bibline,
{ Macdonald's evaluation conjectures and
difference Fourier transform},
Inventiones Math. (1995).

[C4]
\bibline,
{A unification of Knizhnik--Zamolodchikov
and Dunkl operators via affine Hecke algebras},
Inventiones Math. {  106}:2,  (1991), 411--432.

[C5]
\bibline,
{ Integration of quantum many- body problems by affine
Knizhnik--Za\-mo\-lod\-chi\-kov equations},
Pre\-print RIMS--{  776} (1991),
(Advances in Math. (1994)).

[C6]
\bibline,
{Elliptic quantum many-body problem and double affine
 Knizhnik - Zamolodchikov equation},
Commun. Math. Phys. (1995).

[D]
\name{C.F. Dunkl},
{ Hankel transforms associated to finite reflection groups},
Contemp. Math. {138} (1992), 123--138.

[EK]
\name{P.I. Etingof}, and {A.A. Kirillov,Jr.}
{Representation-theoretic proof of the inner product and
symmetry identities for Macdonald's polynomials},
Compositio Mathematica (1995).

[FV]
\name{G. Felder}, and \name{A.P. Veselov}
{Shift operators for the quantum
Calogero - Sutherland problems via Knizhnik - Zamolodchikov
equation}, Communs Math. Phys. {160} (1994), 259--274.

[J]
\name{M.F.E. de Jeu},
{The Dunkl transform }, Invent. Math. {113} (1993), 147--162.

[He]
\name{G.J. Heckman},
{  An elementary approach to the hypergeometric shift operators of
Opdam}, Invent.Math. {  103} (1991), 341--350.
Comp. Math. {  64} (1987), 329--352.
v.Math.{  70} (1988), 156--236.

[HO]
\name{G.J. Heckman}, and \name{E.M. Opdam},
{Yang's system of particles and Hecke algebras},
Preprint (1995).

[K]
\name {V.G. Kac},
{Infinite dimensional Lie algebras},
Cambridge University Press, Cambridge (1990).

[Ki]
\name {A. Kirillov, Jr.},
{Inner product on conformal blocks and Macdonald's
polynomials at roots of unity}, Preprint (1995).

[KL1]
\name{D. Kazhdan}, and \name{ G. Lusztig},
{  Proof of the Deligne-Langlands conjecture for Hecke algebras},
Invent.Math. {  87}(1987), 153--215.

[KL2]
\bibline,
{ Tensor structures arising from affine Lie algebras. III},
J. of AMS {  7}(1994), 335-381.

[KK]
\name{B. Kostant}, and \name{ S. Kumar},
{  T-Equivariant K-theory of generalized flag varieties,}
J. Diff. Geometry{  32}(1990), 549--603.

[M1]
\name{I.G. Macdonald}, {  A new class of symmetric functions },
Publ.I.R.M.A., Strasbourg, Actes 20-e Seminaire Lotharingen,
(1988), 131--171 .

[M2]
\bibline, {  Orthogonal polynomials associated with root
systems},Preprint(1988).

[M3]
\bibline, { Affine Hecke algebras and orthogonal polynomials},
S\'eminaire Bourbaki{  47}:797 (1995), 01--18.

[O1]
\name{E.M. Opdam},
{  Some applications of hypergeometric shift
operators}, Invent.Math.{  98} (1989), 1--18.

[O2]
\bibline, {Harmonic analysis for certain representations of
graded Hecke algebras},
Preprint Math. Inst. Univ. Leiden W93-18 (1993).

\endreferences

\bye